\documentclass{article}
\usepackage{graphicx} 
\usepackage{url}
\usepackage{appendix}

\usepackage{amssymb}
\usepackage{amsthm} 
\usepackage{amsmath} 
\usepackage{booktabs} 
\usepackage{array} 
\usepackage{paralist} 
\usepackage{verbatim} 

\usepackage{subcaption}
\usepackage{xcolor}
\usepackage{multirow}
\usepackage{tikz}
\usetikzlibrary{arrows.meta}
\usepackage[ruled,vlined,linesnumbered]{algorithm2e}
\usepackage{hyperref}

\usepackage{tikz}
\usetikzlibrary{trees}
\usetikzlibrary{positioning}

\usepackage{fancyhdr} 

\title{GPU-Accelerated Orbit Propagation with High-Fidelity Solar Radiation Pressure Modeling}
\author{
Leandro Zardaín$^{1}$,
Ariadna Farrés$^{2}$,
Anna Puig$^{1}$,\\
Kevin Muntal$^{1}$
and Albert Mir$^{1}$
}
\date{August 2026}

\begin{document}

\maketitle

\begin{center}
$^{1}$ Universitat de Barcelona, Barcelona, Spain\\
$^{2}$ University of Maryland, Baltimore County, Baltimore, MD, USA
\end{center}

\begin{abstract}

Solar Radiation Pressure (SRP), the force exerted by photons emitted by the Sun, is one of the main non-gravitational perturbations affecting spacecraft trajectories, making its accurate modeling essential for high-fidelity orbit propagation. While physically-based ray-tracing models improve SRP accuracy, their computational cost becomes a limitation during numerical integration, where the SRP force must be evaluated repeatedly throughout the trajectory propagation.

This paper investigates the integration of high-fidelity SRP models into orbit propagation through two complementary contributions: a Vulkan-based GPU implementation that accelerates direct SRP evaluation, and an extension of the SRP model to account for the dynamic orientation of solar panels. These contributions are evaluated independently through orbit propagation experiments, while a precomputed SRP interpolation strategy (SPAD) is included as an alternative approach for reducing computational cost through offline sampling and interpolation.

Numerical validation shows that the Vulkan implementation preserves the accuracy of  the original OpenGL-based method, showing relative differences below $5\times10^{-4}$  while achieving speed-ups of up to 9.4 for individual SRP computations and up to 15.2 for complete orbit propagation, particularly for geometrically complex spacecraft.

The movable solar panel model shows that neglecting panel motion can produce significant long-term propagation errors, especially for spacecraft with large articulated solar panels, while introducing only a moderate computational overhead.

Based on the experimental evaluation, this work concludes with practical guidelines for integrating high-fidelity SRP models into orbit propagation frameworks, identifying the scenarios under which interpolation-based methods are sufficient and when online high-fidelity SRP computation is justified according to the required balance between computational efficiency, numerical accuracy, and physical fidelity.

\end{abstract}

\section{Introduction}


Recent years have seen a rapid increase in the number and complexity of space exploration missions. Between 2020 and 2030, approximately 150 missions are expected to be launched in Low Earth Orbit, 95 to the Moon and 11 to Mars\footnote{\url{https://www.statista.com/statistics/1169834/space-exploration-missions-worldwide-type/}}. The cost of these missions can reach billions of dollars \footnote{\url{https://www.esa.int/Science_Exploration/Human_and_Robotic_Exploration/International_Space_Station/How_much_does_it_cost}} and therefore, it has become increasingly important to reduce mission costs and improve the efficiency of orbital maneuvers. Therefore, accurate and efficient orbit propagation has become essential for mission design, enabling more precise maneuver planning and reducing propellant consumption. 


Several factors 
affect the orbit determination of a satellite. Among the most important are gravitational perturbations from planets 
and Solar Radiation Pressure (SRP), which can become one of the dominant non-gravitational forces, especially for missions operating close to the Sun.
%
Solar Radiation Pressure is produced when solar photons transfer momentum to the spacecraft surfaces through absorption or reflection. 
This work focuses on missions in which SRP plays a significant role in orbit propagation.

The most common models to approximate the SRP are the Cannonball \cite{Montenbruck} and the N-Plate \cite{Vallado2001}. The Cannonball model assumes the spacecraft is a sphere, whereas the N-Plate model approximates the spacecraft as a collection of flat plates. Their low computational cost makes them attractive for orbit propagation. However, this efficiency is based on a simplified spacecraft geometry and light transport assumptions, simplifying effects such as articulated solar panels, self-occlusions, and multiple reflections. As a result, their accuracy decreases for complex spacecraft geometries and even more so during long-term orbit propagation, where
small errors in the SRP estimation can be accumulated over time, leading to noticeable trajectory deviations and affecting maneuver planning.


In these scenarios, ray tracing is one of the most accurate approaches for high-fidelity SRP computation. It represents the spacecraft as a triangular mesh and it also provides self-shadowing and  multiple reflections. 
Although this approach allows highly accurate SRP estimation, its computational cost becomes a major challenge when integrated into numerical orbit propagation. Unlike standalone SRP evaluations, numerical propagators require thousands or even millions of SRP computations along the trajectory, making execution time a critical factor. Furthermore, modern spacecraft frequently include articulated solar panels whose orientation changes throughout the mission, requiring the spacecraft geometry to be updated before each SRP evaluation.

Our previous work introduced HiFi-SoRaP \cite{simulatorSRP}, a GPU-based ray-tracing framework for high-fidelity SRP computation.  That work focused on accurate force computation and proposed several numerical optimizations, 
including Kahan compensated summation and optimized GPU execution strategies. However, it did not address the challenges associated with repeated SRP evaluation during orbit propagation. 
Moreover, mission analysis tools, such as GMAT \footnote{\url{https://sourceforge.net/projects/gmat/}}, support analytical SRP models for trajectory propagation but do not natively perform on-the-fly high-fidelity ray-tracing SRP evaluation during numerical integration. 
As a result, they do not provide a unified environment for systematically evaluating analytical, interpolation-based, and online high-fidelity SRP models under identical propagation conditions.

This paper extends HiFi-SoRaP into a complete orbit propagation framework. Two complementary strategies are introduced. First, a Vulkan-based implementation substantially accelerates direct online SRP evaluation during numerical integration. Second, movable solar panels are incorporated into the spacecraft model, allowing the SRP computation to account for changes in spacecraft geometry throughout the propagation. To assess the practical applicability of direct high-fidelity SRP computation, the proposed framework is used to compare direct SRP evaluation with SPAD that represents an offline precomputed solution, and provides an efficient baseline to compare interpolation-based and online high-fidelity SRP computation under identical propagation conditions.



In summary, the proposed framework enables a systematic comparison between interpolation-based and online high-fidelity SRP computation within the same orbit propagation environment. Based on this comparison, we derive practical guidelines for selecting the most appropriate SRP strategy according to the required balance between computational efficiency, numerical accuracy and physical fidelity.

\section{Related Work}

Trajectory propagation under perturbative forces has been extensively studied in orbital mechanics and spacecraft dynamics. When non-gravitational perturbations are included, analytical solutions are generally unavailable and the equations of motion must be solved numerically. Consequently, the accuracy of long-term orbit prediction depends on both the numerical integration scheme and the fidelity of the perturbation models employed~\cite{Dormand1980AFO,Hairer2006GNI}.

Among non-gravitational perturbations, Solar Radiation Pressure (SRP) is one of the dominant forces affecting high-altitude Earth satellites, GNSS constellations, and deep-space missions ~\cite{Vallado2001,Montenbruck}.  As mission accuracy requirements have increased, SRP models have evolved from simplified analytical formulations towards high-fidelity approaches capable of representing spacecraft geometry, optical properties, and self-shadowing effects.

For GNSS satellites, Chang et al. \cite{chang} evaluate several SRP models for GPS satellites, including box-wing representations, demonstrating the strong influence of SRP modeling on precise orbit determination. Similarly, in the context of deep-space missions, Kato et al. \cite{rodriguez_rosseta} present a precise modeling of solar and thermal accelerations for the Rosetta mission, highlighting the importance of accurate radiation force estimation for long-term trajectory prediction.

To improve the accuracy of SRP estimation, several authors have proposed ray-tracing approaches that explicitly simulate photon interactions with complex spacecraft geometries. Darugna et al. \cite{li} introduce a ray-tracing SRP modeling approach for the QZS-1 satellite, showing that geometric ray-based methods provide improved accuracy compared to simplified analytical models. Li et al. \cite{zhang} propose a fast SRP modeling technique using ray tracing with multiple reflections and acceleration structures, emphasizing computational efficiency for complex spacecraft geometries.

Although these high-fidelity models improve SRP estimation, they considerably increase the computational cost. 
As ray-surface interactions are independent, SRP computation is easy to parallelize in GPU architectures. Early GPU implementations relied on graphics-oriented APIs such as OpenGL for accelerating visibility analysis and SRP computation \cite{Tichy2014FastFS,HifiSorap,Kenneally16}. More recently, developments incorporated advanced ray-tracing techniques and physically based surface interaction models to improve both accuracy and efficiency \cite{kenneally2018brdf,Kenneally20}. Zardaín et al. \cite{HifiSorap} proposed the HiFi-SoRaP framework, integrating multiple SRP models and GPU-based ray tracing within a common environment for SRP analysis and visualization.

In recent years, modern low-level GPU APIs have been explored to further reduce the computational cost of direct SRP evaluation. Muntal \cite{vulkan} demonstrated that a Vulkan-based ray-tracing implementation accelerates high-fidelity SRP computation by reducing driver overhead and providing explicit GPU resource management.

An alternative strategy to reduce the computational cost of high-fidelity SRP consists of precomputing it over a set of spacecraft attitudes and storing the results in interpolation databases (SPAD). During the orbit propagation, SRP forces are then obtained through interpolation instead of executing the ray-tracing algorithm at every integration step. These approaches significantly reduce the computational cost while introducing an approximation whose accuracy depends on the sampling strategy and interpolation scheme.

Besides computational efficiency, another challenge for high-fidelity SRP computation is accurately representing spacecraft configurations that evolve during the mission. Mir \cite{mobilePanels} proposed a method to model articulated solar panels by updating the spacecraft geometry on the CPU before each GPU ray-tracing evaluation. Their work addressed the geometric reconfiguration required for articulated spacecraft panels but did not consider long-term numerical orbit propagation or the integration of SRP evaluation within a trajectory propagation framework.

In summary, existing work has therefore addressed either efficient high-fidelity SRP computation or improved geometric modeling, but these developments have rarely been evaluated together within a complete orbit propagation framework. Furthermore, most comparisons focus on isolated SRP computations rather than on their cumulative impact during long-term numerical propagation. In this paper, we investigate three complementary strategies for integrating high-fidelity SRP models into orbit propagation: GPU acceleration through a Vulkan-based implementation, improved physical realism through movable solar panel modeling, and an interpolation-based approach (SPAD) for reducing the cost of repeated SRP evaluations.

\section{Numerical Challenges of High-Fidelity Orbit Propagation}
\label{sec:framework}

This section introduces the mathematical formulation of the orbit propagation problem and discusses the main challenges associated with incorporating high-fidelity Solar Radiation Pressure (SRP) models into numerical orbit propagation. These challenges motivate the research questions addressed in the remainder of the paper.

Specifically, the objective of the orbit propagation model is to determine the future trajectory of a spacecraft from a given initial state by numerically integrating the equations of motion. In this work, the spacecraft dynamics are modeled using the classical two-body problem, augmented with Solar Radiation Pressure (SRP), which is treated as the only perturbing force. The total acceleration is therefore expressed as

\begin{equation}
\mathbf{a}(t)=\mathbf{a}_g(t)+\mathbf{a}_{SRP}(t).
\end{equation}

\noindent where $\mathbf{a}_g$ denotes the gravitational acceleration and $\mathbf{a}_{SRP}$ represents the acceleration generated by the selected SRP model. Moreover, the spacecraft equations of motion can be expressed in terms of the state vector
\begin{equation}
\mathbf{y}(t)=
\begin{bmatrix}
\mathbf{x}(t)\\
\mathbf{v}(t)
\end{bmatrix},
\end{equation}

\noindent where $\mathbf{x}(t)$ and $\mathbf{v}(t)$ denote the position and velocity vectors, respectively. The equations of motion are then written as

\begin{equation}
\dot{\mathbf{y}}(t)=
\begin{bmatrix}
\mathbf{v}(t)\\
\mathbf{a}(t)
\end{bmatrix}
=\mathbf{f}(t,\mathbf{y}).
\end{equation}

Since no analytical solution exists once SRP is included, the spacecraft trajectory must be obtained through the numerical integration of this system of differential equations. 

The spacecraft trajectory is propagated using an explicit Runge-Kutta 7/8 (RK78) integration scheme. RK78 \cite{Dormand1980AFO} evaluates the state derivative multiple times within each integration step to achieve high numerical accuracy. Consequently, every integration step requires multiple evaluations of the SRP acceleration.

For analytical SRP models, such as the Cannonball or N-plate models, the computational cost of each evaluation is negligible compared with the integration process itself. However, when the SRP acceleration is computed using high-fidelity ray tracing models, every evaluation requires performing a complete ray-tracing simulation over the spacecraft geometry. As a result, the repeated SRP evaluations required by the numerical integrator become the dominant computational cost of the orbit propagation.

Besides computational efficiency, the accuracy of the SRP model also influences the propagated trajectory. Many high-fidelity implementations assume a static spacecraft geometry throughout the propagation. However, for spacecraft with movable solar panels, the panel orientation changes continuously to maintain Sun-pointing conditions, modifying the illuminated surfaces and the resulting radiation forces. Neglecting these geometric changes may therefore introduce significant errors during long-term orbit propagation.

The previous discussion identifies two fundamental challenges for incorporating high-fidelity SRP models into orbit propagation: reducing the computational cost associated with repeated SRP evaluations and improving the physical realism of the spacecraft model. The computational challenge is addressed from two complementary perspectives: replacing direct SRP evaluations with interpolation-based approximations and accelerating direct evaluations through a Vulkan-based implementation. Together with the incorporation of movable solar panels, these objectives motivate the following research questions.

\begin{itemize}

\item \textbf{RQ1.} Can an interpolation-based representation of the SRP replace direct GPU ray-tracing evaluation during orbit propagation while maintaining sufficient propagation accuracy?

\item \textbf{RQ2.} To what extent can a Vulkan-based GPU implementation reduce the computational cost of direct high-fidelity SRP evaluation while preserving numerical accuracy?

\item \textbf{RQ3.} What is the impact of incorporating movable solar panels into the spacecraft model on orbit propagation accuracy and computational cost?
\end{itemize}

The following section presents the proposed orbit propagation framework developed to address these research questions.

\section{Proposed Orbit Propagation Framework}

This section presents the proposed orbit propagation framework
developed to address the research questions introduced in the previous section.
The framework extends the existing HiFi-SoRaP framework.
HiFi-SoRaP evaluates Solar Radiation Pressure (SRP) forces for individual spacecraft configurations and supports the comparison of analytical models, such as the Cannonball and N-Plate models, with high-fidelity GPU ray-tracing models.
Unlike HiFi-SoRaP, which evaluates the SRP acceleration for a single spacecraft state, the proposed framework incorporates these SRP evaluations into a numerical orbit propagator, enabling the comparison of different SRP computation strategies under identical propagation conditions.

\subsection{Workflow Overview}

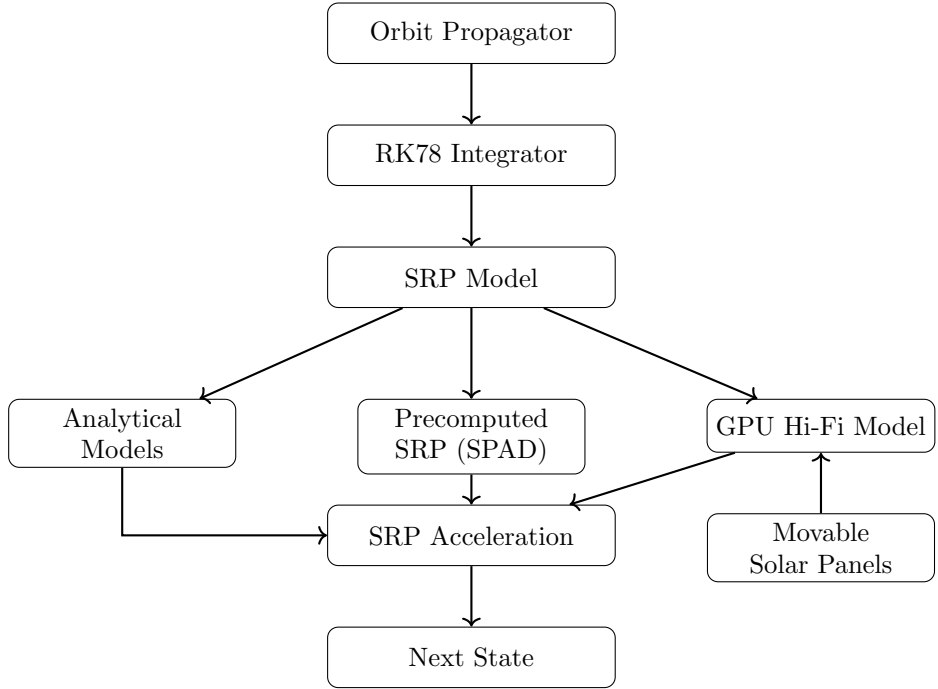
\begin{figure}[ht]
\centering
\begin{tikzpicture}[
    node distance=8mm,
    box/.style={draw, rounded corners, align=center, minimum width=3.8cm, minimum height=8mm},
    smallbox/.style={draw, rounded corners, align=center, minimum width=3.0cm, minimum height=7mm},
    arrow/.style={->, thick}
]

\node[box] (prop) {Orbit Propagator};
\node[box, below=of prop] (rk) {RK78 Integrator};
\node[box, below=of rk] (srp) {SRP Model};

\node[smallbox, below left=12mm and 12mm of srp] (ana) {Analytical\\Models};

\node[smallbox, below=12mm of srp] (spad) {Precomputed\\SRP (SPAD)};

\node[smallbox, below right=12mm and 12mm of srp] (rt) {GPU Hi-Fi Model};

\node[smallbox, below=of rt] (mp) {Movable\\Solar Panels};

\node[box, below=26mm of srp] (acc) {SRP Acceleration};

\node[box, below=of acc] (next) {Next State};

\draw[arrow] (prop) -- (rk);
\draw[arrow] (rk) -- (srp);

\draw[arrow] (srp) -- (ana);
\draw[arrow] (ana) |- (acc);

\draw[arrow] (srp) -- (spad);
\draw[arrow] (spad) -- (acc);

\draw[arrow] (srp) -- (rt);
\draw[arrow] (rt) -- (acc);

\draw[arrow] (mp) -- (rt);

\draw[arrow] (acc) -- (next);

\end{tikzpicture}
\caption{Workflow of the proposed orbit propagation framework.}
\label{fig:framework}
\end{figure}

Figure~\ref{fig:framework} illustrates the overall workflow of the proposed framework. The orbit propagator employs an RK78 numerical integrator, which repeatedly requests the SRP acceleration during each integration step. Depending on the selected model, the SRP acceleration can be computed either using analytical models, the SPAD interpolation model, or by means of the GPU ray-tracing framework. When movable solar panels are enabled, the spacecraft geometry is first updated according to the panel orientation before the standard ray-tracing algorithm is executed. The resulting SRP acceleration is then returned to the numerical integrator to advance the spacecraft state.

This modular organization allows different SRP computation methods to be exchanged without modifying the numerical integration process, enabling direct comparisons under identical propagation conditions.
All SRP models share the same numerical integration scheme and propagation settings, ensuring that any observed differences are exclusively attributable to the SRP computation strategy.

\subsection{Interpolation-Based SRP Computation}

To reduce the computational cost associated with repeated SRP evaluations, the framework also supports the SPAD interpolation model. Instead of executing a complete GPU ray-tracing simulation during every integration step, the SRP acceleration is obtained from a precomputed interpolation dataset generated offline. This allows orbit propagation to be performed with a computational cost comparable to analytical models while preserving much of the accuracy of the original high-fidelity model.

\subsection{Vulkan-based SRP Computation} \label{sec:vulkan}

The proposed framework incorporates a Vulkan-based GPU backend for high-fidelity SRP evaluation. The original OpenGL-based implementation was adapted to a Vulkan compute pipeline~\cite{vulkan}. The proposed Vulkan backend preserves the original SRP formulation while replacing the OpenGL compute pipeline with a Vulkan implementation that provides explicit control over memory management, synchronization, and compute scheduling. These capabilities reduce the overhead associated with repeated SRP evaluations during orbit propagation. The numerical improvements introduced in the HiFi-SoRaP framework \cite{simulatorSRP}, including Kahan compensated summation, were incorporated into this implementation to maintain numerical consistency with the rest of the framework.

The Vulkan version preserves the structure of the original ray tracing algorithm but introduces several optimizations. The computation is still performed per pixel of the discretized image plane, where each pixel corresponds to an independent ray. This structure allows a fully parallel execution of the SRP computation on the GPU. Vulkan provides explicit control over memory transfers, synchronization, and work scheduling, allowing the repeated SRP evaluations required during orbit propagation to be executed more efficiently.
Each pixel is processed independently through compute shaders, and intermediate results are stored in GPU memory.

Secondary rays are handled iteratively following the same physical assumptions as in the original model. The implementation also supports different levels of scattering depending on the number of secondary rays generated at each interaction.

The final SRP acceleration is then obtained by aggregating all contributions. To reduce computational overhead, a parallel reduction strategy is used to sum the contributions of all pixels directly on the GPU, avoiding unnecessary transfers to the CPU. This is achieved using subgroup operations, which allow partial sums to be computed efficiently within each execution group before a final aggregation step.

\subsection{SRP Computation with Movable Solar Panels}\label{sec:PM}

One of the main limitations of the original HiFi-SoRaP framework was the assumption of a rigid spacecraft geometry throughout the SRP computation. While this assumption is adequate for many spacecraft, satellites equipped with steerable solar panels require their geometry to be updated continuously according to the Sun direction.

To overcome this limitation, Mir \cite{mobilePanels} extended the spacecraft model to support articulated solar panels. The numerical improvements introduced in the HiFi-SoRaP framework \cite{simulatorSRP}, including Kahan compensated summation, were incorporated into this implementation to ensure numerical consistency with the rest of the framework. The graph representation required for the proposed approach is summarized below. A complete mathematical formalization is available in \cite{mobilePanels}.


\subsubsection{Tree Definition}

To support articulated structures, the spacecraft is modeled as a rooted directed tree 
\[
G=(N,E,L),
\]

\noindent where \(N\) is the set of nodes, \(E\) the set of directed edges defining the kinematic hierarchy, and \(L\) the set of edge labels describing the rotational constraints associated with each joint. Within this representation, each node corresponds to a rigid group of spacecraft components, while edges describe the kinematic relationships between them.
Formally, each node is defined as
\[
n=(O_n,t_n),
\]

\noindent where \(O_n\) is the subset of spacecraft mesh parts associated with the node, and \(t_n\) is a binary flag indicating whether the node corresponds to an articulated solar panel. The sets of mesh parts form a partition of the spacecraft model, i.e., every mesh part belongs to exactly one node. Each edge $
e=(n_p,n_q)\in E$ represents a parent-child relationship in the kinematic hierarchy. 

Since a rotation applied to an articulated component also affects all components attached to it, the notion of a subtree is required. For a node \(n\in N\), the subtree rooted at \(n\) contains the node itself together with all its descendants. Any rotation applied through the edge connecting \(n\) to its parent affects every mesh component belonging to this subtree.

This hierarchical representation allows rigid transformations to be applied efficiently to complete articulated structures by operating on subtrees instead of individual mesh elements. The rotational properties associated with each joint are stored as edge labels. Each edge \(e\in E\) is associated with a label $l=(p,A,I)$, where \(p\) is the pivot point, \(A=\{\mathbf{a}_1,\ldots,\mathbf{a}_k\}\) is the set of rotation axes, and \(I=\{I_1,\ldots,I_k\}\) is the corresponding set of admissible angular intervals.

\subsubsection{SRP Computation}

Once the spacecraft has been represented using this hierarchical model, the geometry can be updated dynamically before each SRP evaluation.

During each computation, the algorithm traverses the tree and identifies movable joints. For each of these, the following steps are applied:

\begin{enumerate}
    \item The subtree starting at the corresponding node is extracted.
    \item The normal of the solar panels in that subtree is computed.
    \item The panels are rotated so that their normal is aligned as closely as possible with the Sun direction.
    \item The rotation is clamped according to the predefined angular limits.
    \item The transformation is applied to all triangles in that subtree before computing SRP.
\end{enumerate}

This enables the solar panels to continuously adjust their orientation towards the Sun direction at each simulation step, improving the physical realism of the simulation.

After updating the spacecraft geometry, the standard ray-tracing algorithm is executed without further modifications. The numerical improvements incorporated into the HiFi-SoRaP framework~\cite{simulatorSRP}, including Kahan compensated summation, are also used in this implementation to ensure numerical consistency across all SRP computation strategies.
Consequently, the proposed extension modifies only the spacecraft geometry prior to SRP evaluation, while preserving the original ray-tracing algorithm.

\section{Experimental Setup}

To ensure a fair comparison, all implementations are assessed under identical simulation conditions, allowing any differences in measured performance to be attributed exclusively to the SRP computation strategy.

The experiments were executed on an AMD Ryzen 5 3600X processor with 16 GB RAM and an NVIDIA GeForce RTX 2060 Super GPU.

\subsection{Spacecraft Models and Mission Definition}

Two spacecraft models are considered to evaluate the proposed framework: a GPS-inspired BoxWing satellite and a high-fidelity model of the Parker Solar Probe (PSP). These models represent two different levels of geometric complexity.

Table~\ref{table:datasetsVulkanPM} summarizes the main properties of both spacecraft.

\begin{table}[htbp]
\centering
\caption{Summary of the spacecraft models used in the experiments.}
\label{table:datasetsVulkanPM}
\small
\begin{tabular}{lccc}
\toprule
\textbf{Spacecraft} & \textbf{Geometry} & \textbf{Materials} & \textbf{Dry Mass} \\
\midrule

\begin{tabular}[c]{@{}c@{}} 
PSP \\
\includegraphics[width=0.25\textwidth]{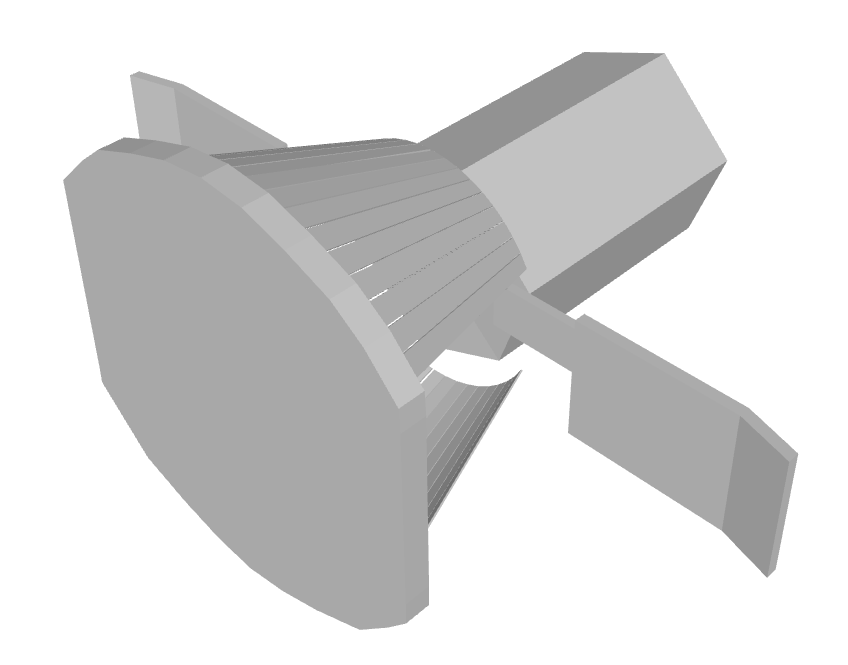}
\end{tabular}
&
\begin{tabular}[c]{@{}l@{}}
Num. Triangles = 39257 \\
Size: 3.7 $\times$ 3.4 $\times$ 4 (m)
\end{tabular}
&
\begin{tabular}[c]{@{}l@{}}
\textbf{solar shield} \\
$p_s = 0.3$ \\
$p_d = 0.6$ \\
$p_a = 0.1$ \\
\textbf{bus} \\
$p_s = 0.3$ \\
$p_d = 0.6$ \\
$p_a = 0.1$ \\
\textbf{solar array} \\
$p_s = 0.176$ \\
$p_d = 0.044$ \\
$p_a = 0.78$
\end{tabular}
&
555 kg \\
\midrule

\begin{tabular}[c]{@{}c@{}}
BoxWing \\
\includegraphics[width=0.22\textwidth]{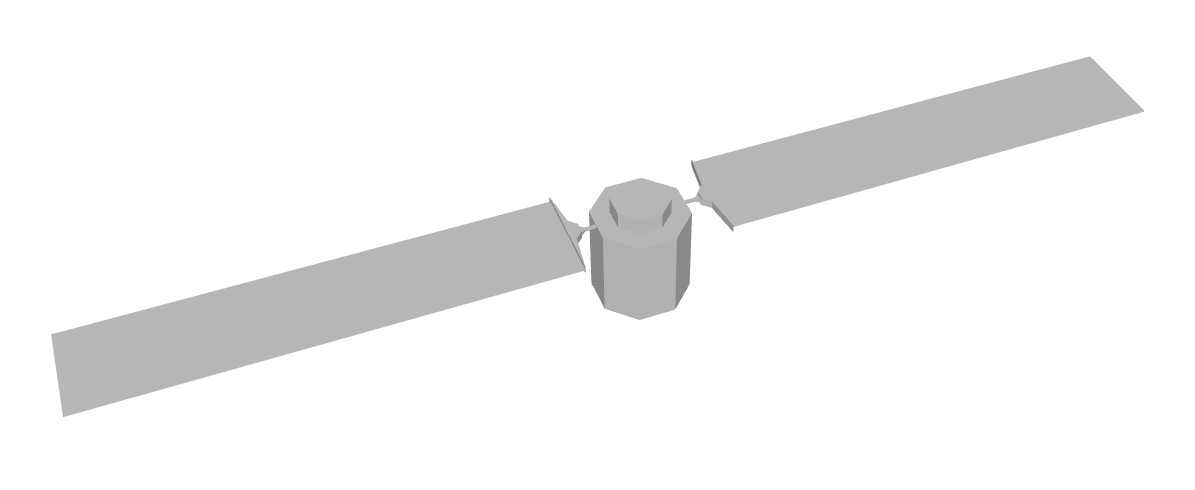}
\end{tabular}
&
\begin{tabular}[c]{@{}l@{}}
Num. Triangles = 212 \\
Size: 21.9 $\times$ 2.36 $\times$ 2.07 (m)
\end{tabular}
&
\begin{tabular}[c]{@{}l@{}}
\textbf{bus} \\
$p_s = 0.4$ \\
$p_d = 0.1$ \\
$p_a = 0.5$ \\
\textbf{solar array} \\
$p_s = 0.05$ \\
$p_d = 0.20$ \\
$p_a = 0.75$
\end{tabular}
&
2269 kg \\

\bottomrule
\end{tabular}
\end{table}

All experiments use the same heliocentric elliptical orbit in order to isolate the computational effects introduced by the different SRP implementations. The orbit is defined by a semi-major axis of 0.557 AU, an eccentricity of 0.682, and an inclination of $0^\circ$. The right ascension of the ascending node (RAAN) and the argument of periapsis are also fixed at $0^\circ$. Throughout the simulation, the spacecraft attitude remains fixed, eliminating attitude changes as a source of variability.

\begin{figure}[ht!]

\includegraphics[width=0.8\textwidth]{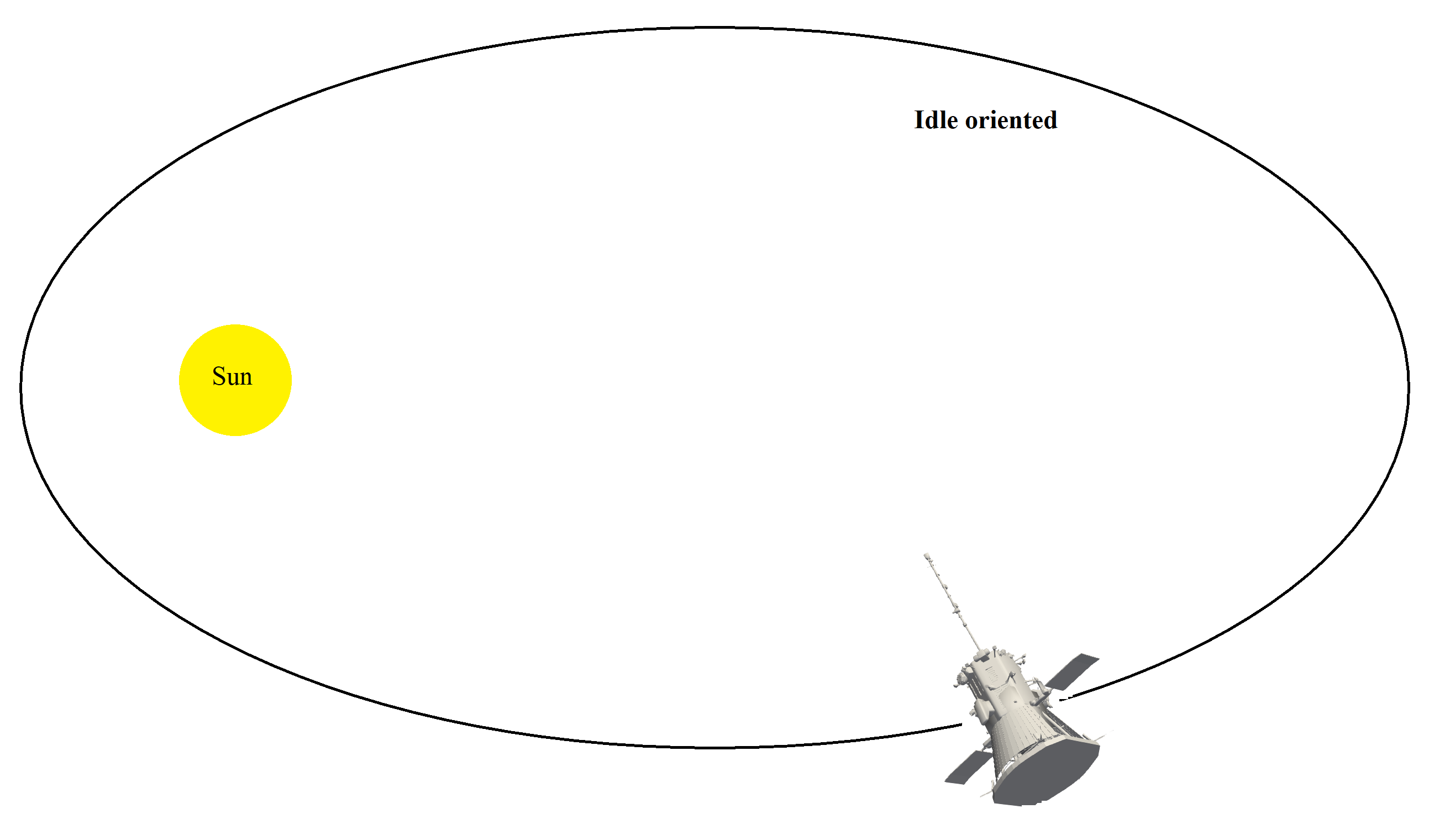}
\caption{Orbit defined for this experiment.}
\label{fig:simulatorsrp-orbital}
\end{figure}

Figure~\ref{fig:simulatorsrp-orbital} illustrates the orbital configuration adopted throughout this section.

\subsection{Evaluation Metrics}

The proposed framework is evaluated using accuracy and performance metrics. Accuracy metrics quantify differences in the computed SRP forces and propagated trajectories, whereas performance metrics measure the computational cost of SRP evaluations and complete orbit propagation.

\begin{itemize}
    \item Accuracy metrics ($A_i$), which quantify differences in SRP forces and propagated trajectories.
    \item Performance metrics ($P_i$), which measure the computational cost associated with SRP evaluations and complete trajectory simulations.
\end{itemize}

The evaluated metrics are
\begin{itemize}
    \item $A_1$: Mean squared error (MSE), maximum error, and relative mean squared error (MSER) of the computed SRP forces over all SRP evaluations.
    
    \item $A_2$: Euclidean distance between propagated trajectories, reported as minimum, maximum, and mean separation throughout the simulation.

    \item $P_1$: Computational cost of SRP force evaluations, reported as total and mean execution times.
    
    \item $P_2$: Total computational cost required to simulate the complete trajectory.
\end{itemize}

\section{Experimental Methodology}

Three complementary experiments are performed to evaluate both proposed strategies.

For the SRP evaluation metrics ($A_1$ and $P_1$), the SRP acceleration is computed for 1,296 uniformly distributed Sun-spacecraft directions covering the entire sphere. This experiment evaluates the numerical accuracy and computational cost of the SRP computation independently of orbit propagation.

For the trajectory metrics ($A_2$ and $P_2$), the spacecraft is propagated along the reference heliocentric elliptical orbit described in the previous subsection for 20 complete orbital revolutions. The same initial conditions, numerical integration scheme, and propagation parameters are used for all SRP models.

\subsection{RayTrace GPU vs SPAD Strategy}


The objective of this study is to evaluate the trade-off between accuracy and computational performance when using different SRP modeling approaches for orbit propagation. In particular, the analysis compares direct SRP computation using RayTrace GPU with precomputed SPAD datasets through interpolation to determine whether the performance improvements achieved through interpolation justify the approximation error introduced with respect to direct high-fidelity SRP evaluation.

The Parker Solar Probe mission described in the previous section is used as the reference case for all simulations. The spacecraft is propagated along the orbital scenarios defined in Table~\ref{table:orbital_parameters} using the following SRP models:

\begin{itemize}
\item No SRP: the SRP acceleration is neglected ($\mathbf{a}_{SRP}=\mathbf{0}$).
\item N-Plate.
\item RayTrace GPU: using a projection window of 256$\times$256 cells, two levels of secondary reflections, and 16 diffuse rays. 
\item SPAD interpolation using two different SPAD datasets. Both datasets were generated using RayTrace GPU with the same configuration described above and a sampling grid of 20 azimuth values and 10 elevation values. The first dataset starts at $-180^{\circ}$ azimuth and $-90^{\circ}$ elevation, whereas the second dataset is shifted by $5^{\circ}$, starting at $-175^{\circ}$ azimuth and $-85^{\circ}$ elevation. These datasets are denoted as SPAD 20x10 and SPAD 20x10S (Shifted), respectively.
\end{itemize}

\begin{table}[htbp]
\centering
\caption{Summary table of all orbits the satellite will follow where C denotes a circular orbit, and E(S), E(M), and E(L) denote small-, medium-, and large-elliptical orbits, respectively.}
\label{table:orbital_parameters}
\footnotesize
\makebox[0.5\textwidth][c]{%
\begin{tabular}{lcccc}
\toprule
\textbf{Type of orbit} & \textbf{C} &
\textbf{E(S)} &
\textbf{E(M)} &
\textbf{E(L)} \\
\midrule
Semi-Major Axis & 0.388 AU & 0.388 AU & 0.437 AU & 0.557 AU \\
Eccentricity & 0 & 0.881 & 0.786 & 0.682 \\
Inclination & $0^\circ$ & $0^\circ$ & $0^\circ$ & $0^\circ$ \\
Right Ascension \\of Ascending Node & $0^\circ$ & $0^\circ$ & $0^\circ$ & $0^\circ$ \\
Argument of periapsis & $0^\circ$ & $0^\circ$ & $0^\circ$ & $0^\circ$ \\
\bottomrule
\end{tabular}
}
\end{table}

The SPAD-based models are used to quantify the error introduced when replacing direct ray tracing evaluations by interpolated SRP accelerations. The SPAD 20x10 dataset contains samples defined at an elevation angle of $0^{\circ}$, which allows the interpolator to perform a linear interpolation along the azimuth direction for the scenarios considered. In contrast, the SPAD 20x10S dataset does not contain samples exactly at $0^{\circ}$ elevation, forcing the interpolator to perform a bilinear interpolation. 

The test cases considered correspond to the orbital and attitude scenarios summarized in Table~\ref{table:all_orbital_simulations}. These include:

\begin{itemize}
\item One circular orbit with fixed spacecraft attitude ($S_0$).
\item Nine elliptical orbit scenarios combining three orbital configurations ($E(S)$, $E(M)$ and $E(L)$) with three attitude transition angles ($S_1$, $S_2$ and $S_3$).
\item One elliptical orbit with continuous attitude variation ($S_4$).
\end{itemize}

\begin{table}[htbp]
\centering
\caption{Summary table of all test case scenarios considered for different orbits and simulations.}
\label{table:all_orbital_simulations}
\small
\begin{tabular}{lcccc}
\toprule
 & \textbf{C} & \textbf{E(S)} & \textbf{E(M)} & \textbf{E(L)} \\
\midrule
$S_{0}$ & X &  &  &  \\
$S_{1}$ &  & X & X & X \\
$S_{2}$ &  & X & X & X \\
$S_{3}$ &  & X & X & X \\
$S_{4}$ &  & X &  &  \\
\bottomrule
\end{tabular}
\end{table}

For each simulation, the resulting trajectories are compared by analyzing both the spacecraft position and the evolution of the corresponding Keplerian orbital elements.

Based on the characteristics of the selected scenarios, two trends are expected regarding the interpolation accuracy. First, for a given attitude-control strategy $S_i$, the largest discrepancies between the SPAD interpolation and the direct RayTrace GPU evaluation are expected for the $E(S)$ orbit. Since the spacecraft remains closer to the Sun in this configuration, the SRP acceleration has a greater influence on the trajectory, making interpolation errors more likely to accumulate into observable differences in the orbital evolution.

Second, the continuous-attitude scenario $S_4$ is expected to produce larger interpolation errors than the discrete-attitude scenarios $S_1$, $S_2$, and $S_3$. During the continuous attitude transition, the SRP acceleration spans a wider range of intermediate illumination conditions, increasing the likelihood of deviations between the interpolated SRP acceleration and the direct RayTrace GPU solution. Consequently, this scenario represents the most demanding case for assessing the accuracy of the interpolation approach.

\subsection{Vulkan Strategy}

The Vulkan strategy compares the original OpenGL-based GPU ray-tracing implementation with the proposed Vulkan-based backend. Both implementations preserve the same SRP formulation, spacecraft models, numerical integration scheme, and propagation parameters. Therefore, any observed differences in the accuracy or performance metrics can be attributed exclusively to the GPU implementation.

\subsection{Movable Solar Panels Strategy}

The movable solar panels strategy compares two spacecraft configurations: one with fixed solar panels and another with movable solar panels that continuously orient toward the Sun during propagation. Both configurations use the same SRP formulation, numerical integration scheme, and propagation parameters. Consequently, any observed differences in the evaluation metrics can be attributed exclusively to the incorporation of movable solar panels.

\label{Simulator}

\section{ Direct GPU Evaluation Versus SPAD Interpolation}

The first research question investigates whether an interpolation-based representation of the SRP acceleration can replace direct RayTrace GPU evaluation during orbit propagation while maintaining sufficient propagation accuracy.

To answer this research question, three aspects are evaluated. First, the propagation accuracy achieved by the interpolation model must be compared with that obtained from direct RayTrace GPU evaluation in order to quantify the interpolation error introduced during orbit propagation. In particular, it is necessary to determine whether this error depends on the orbital scenario or on the spacecraft attitude, since both factors influence the variability of the SRP acceleration.

Second, the computational performance of both approaches must be assessed. Although interpolation is expected to considerably reduce the computational cost of SRP evaluation, this reduction is only meaningful if the associated loss of accuracy remains acceptable for orbit propagation.

Finally, the influence of spacecraft attitude must also be analyzed. Changes in attitude modify the illumination conditions and may introduce rapid variations in the SRP acceleration, potentially affecting the accuracy of the interpolated representation.

The following experiments compare different orbital scenarios, spacecraft attitude configurations, and SRP computation methods. These experiments evaluate the accuracy and computational efficiency of direct RayTrace GPU evaluation and SPAD interpolation under representative operating conditions.

\subsection{Results}

This section compares SPAD-based interpolation with orbit propagation using direct RayTrace GPU evaluation in terms of trajectory accuracy and computational performance.

Eleven test scenarios, including circular and elliptical heliocentric orbits with different spacecraft attitude profiles, were analyzed. Representative results are presented below.

\paragraph{Circular Orbit ($S_0$ - C).}

\begin{figure}[htbp]
\centering
\includegraphics[width=1.\textwidth]{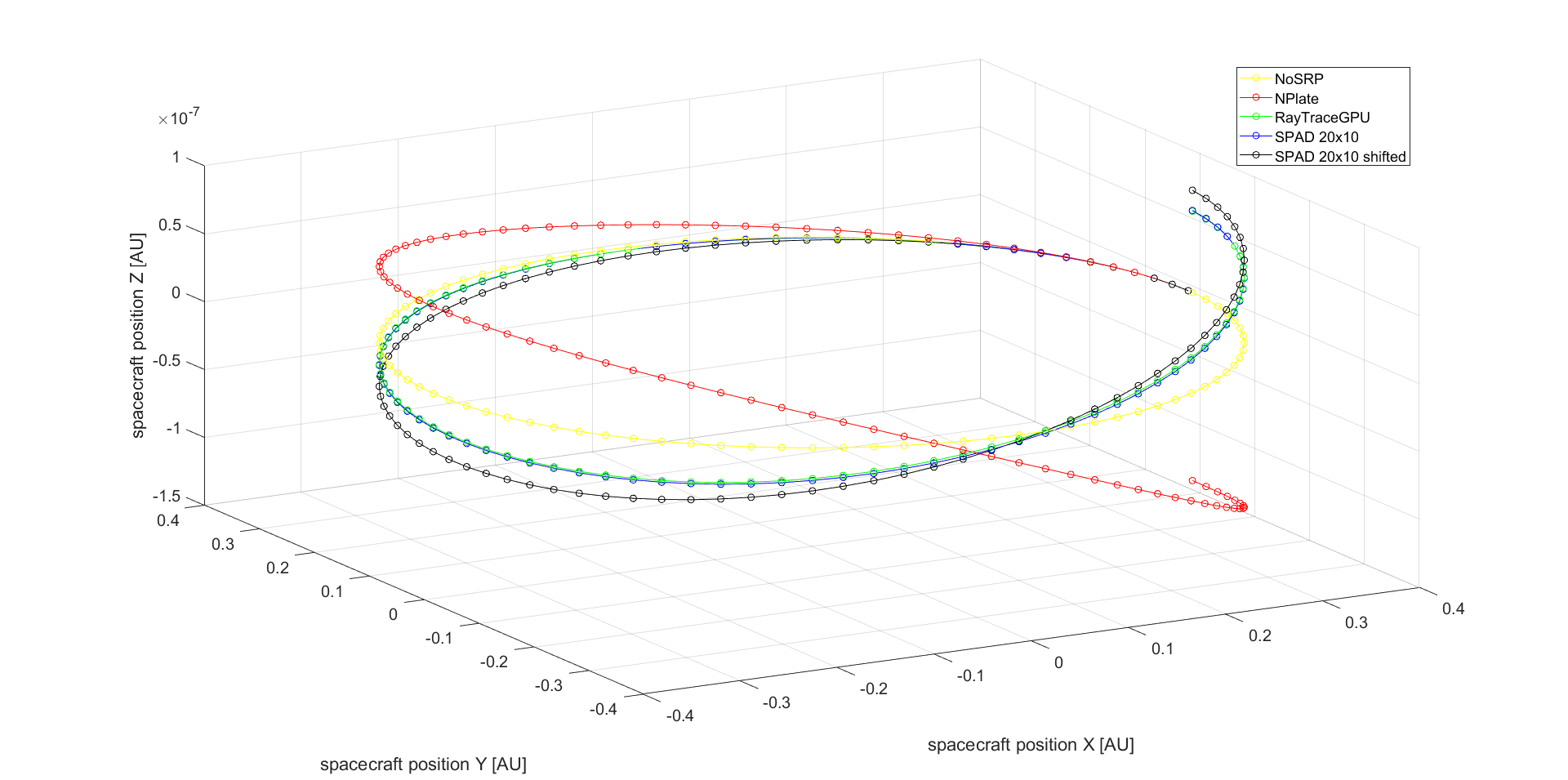}
\caption{Three-dimensional trajectory for the circular-orbit scenario.}
\label{fig:circular_orbit_3d}
\end{figure}

\begin{figure}[htbp]
\centering
\includegraphics[width=1.\textwidth]{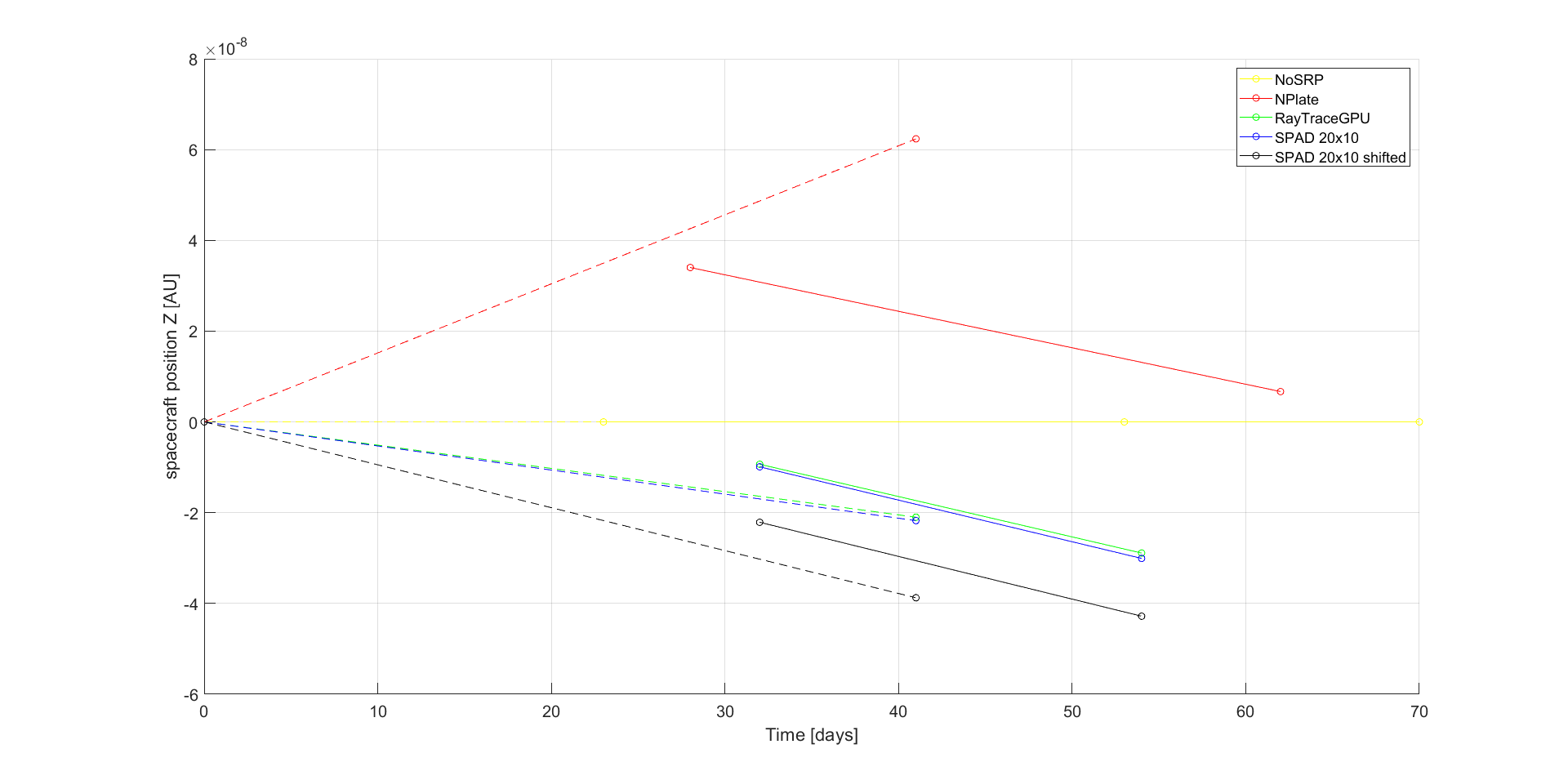}
\caption{Evolution of the spacecraft $z$-coordinate for the circular-orbit scenario.}
\label{fig:circular_orbit_z}
\end{figure}

\begin{figure}[htbp]
\centering
\includegraphics[width=1.\textwidth]{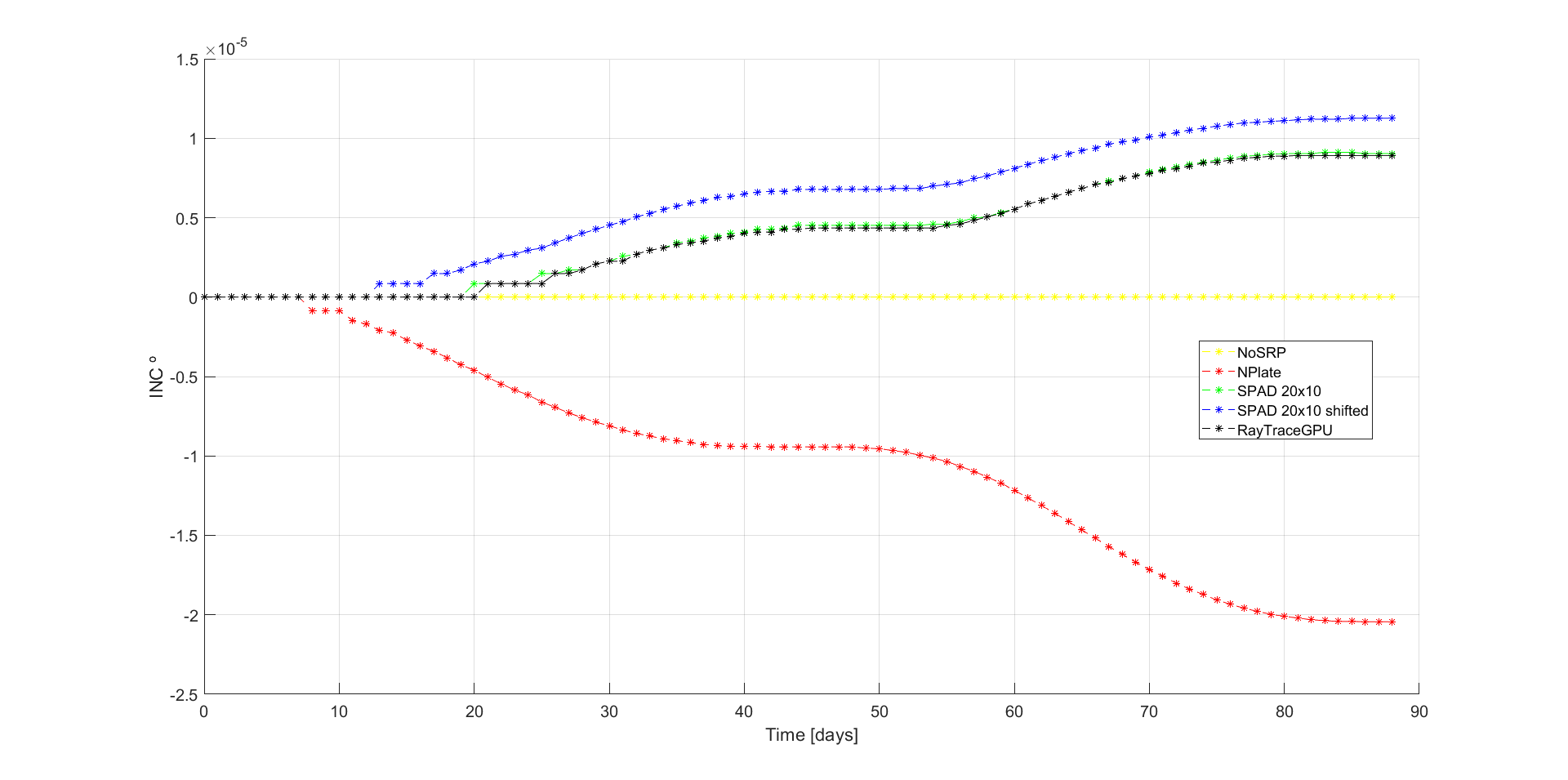}
\caption{Evolution of the orbital inclination for the circular-orbit scenario.}
\label{fig:circular_orbit_inclination}
\end{figure}

Figure~\ref{fig:circular_orbit_3d}-Figure~\ref{fig:circular_orbit_inclination} present the results obtained for the circular-orbit scenario.

Figures~\ref{fig:circular_orbit_3d} and~\ref{fig:circular_orbit_z} show the differences in position between the trajectories generated using the different SRP models. Figure~\ref{fig:circular_orbit_3d} presents the three-dimensional trajectories, while Figure~\ref{fig:circular_orbit_z} shows the evolution of the spacecraft $z$-coordinate at aphelion and perihelion.

The largest deviations appear at the aphelion points. Considering RayTrace GPU as the reference solution, the final position differences are approximately 0.12 km for SPAD 20x10, 2.54 km for SPAD 20x10S, and 12.06 km for N-Plate.

Figure~\ref{fig:circular_orbit_inclination} shows the evolution of the orbital inclination. Although the nominal orbit is planar, SRP introduces small inclination variations. The inclination differences with respect to RayTrace GPU are approximately $1.0\times10^{-7}\,^\circ$  for SPAD 20x10, $1.2\times10^{-6}\,^\circ$  for SPAD 20x10S, and $2.93\times10^{-5}\,^\circ$  for N-Plate.

These results indicate that SPAD 20x10 reproduces the RayTrace GPU trajectory with very high accuracy in this scenario.

\paragraph{Elliptical Orbit Family ($S_1$-$S_3$).}

The family of elliptical orbits exhibits the same general behavior observed in the circular case, although the magnitude of the deviations depends on the orbital size and the attitude profile.

Among all the analyzed configurations, the largest differences are observed for the smallest elliptical orbit, E(S), which corresponds to the trajectory with the closest approach to the Sun. This behavior is consistent with the stronger SRP perturbation experienced near the Sun.

\begin{figure}[htbp]
\centering
\includegraphics[width=1.\textwidth]{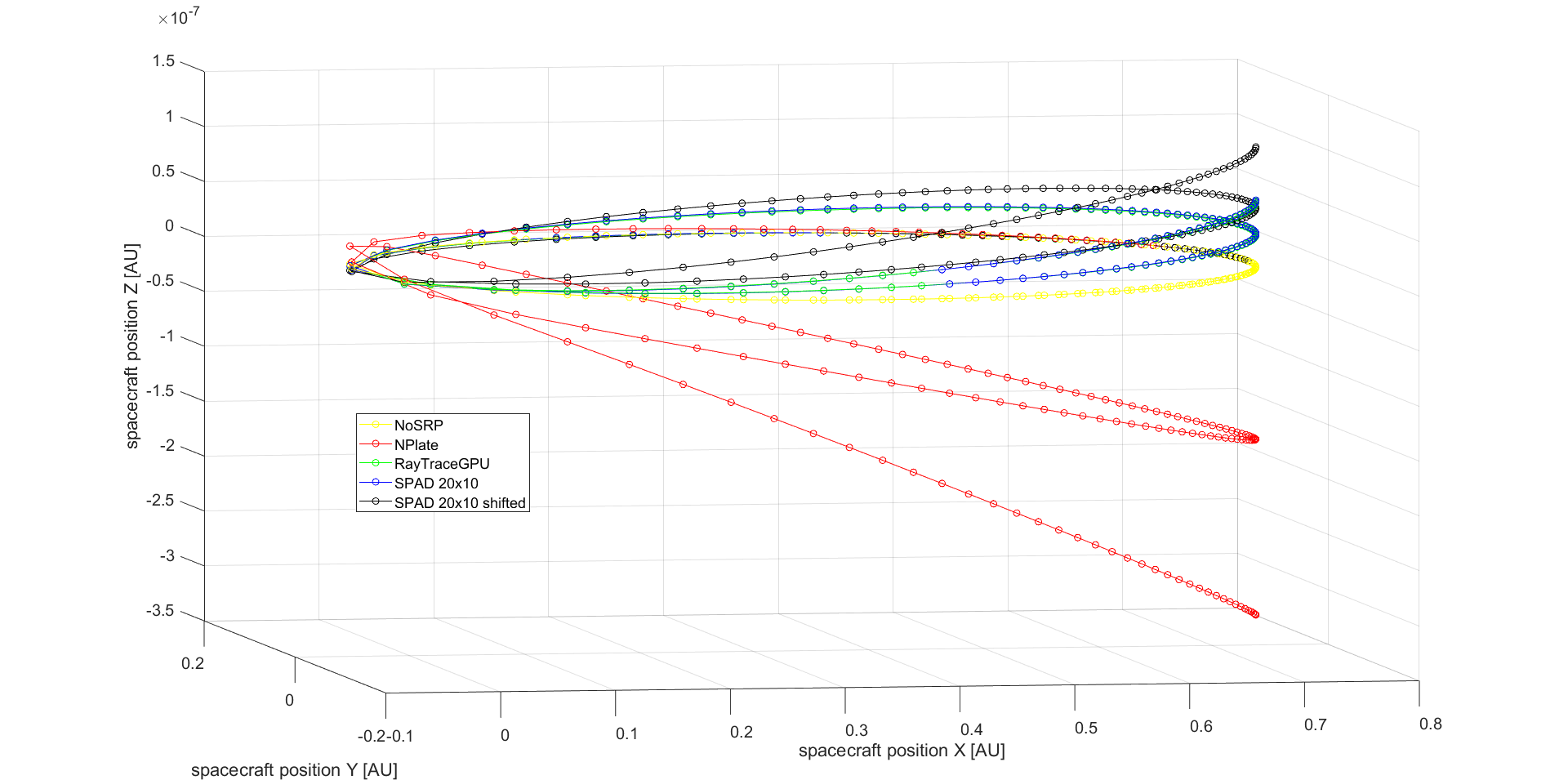}
\caption{Three-dimensional trajectories for scenario $S_2$-E(S).}
\label{fig:closest_normal_3d}
\end{figure}

\begin{figure}[htbp]
\centering
\includegraphics[width=1.\textwidth]{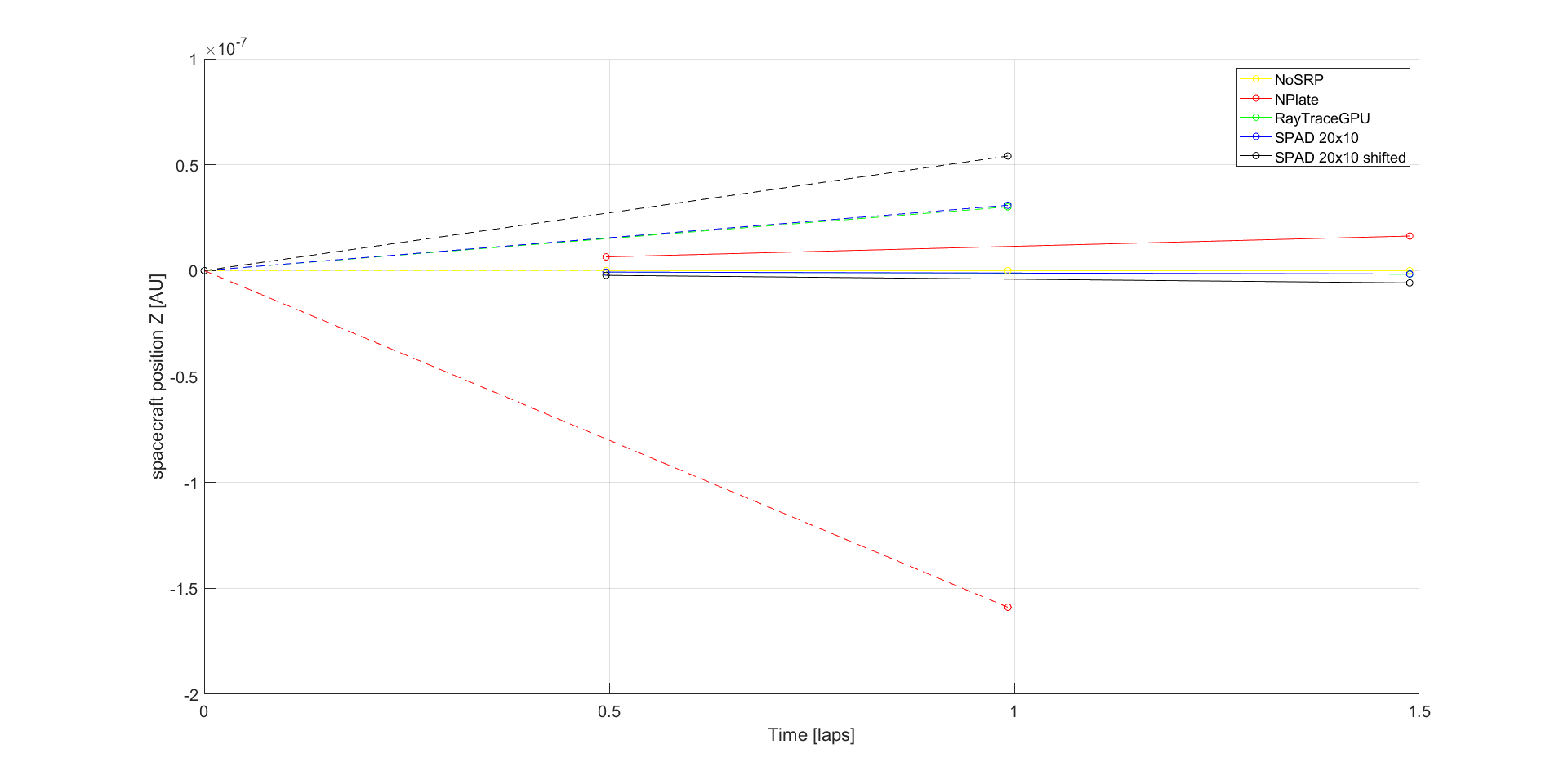}
\caption{Evolution of the spacecraft $z$-coordinate at aphelion and perihelion for scenario $S_2$-E(S).}
\label{fig:closest_normal_z}
\end{figure}

\begin{figure}[htbp]
\centering
\includegraphics[width=1.\textwidth]{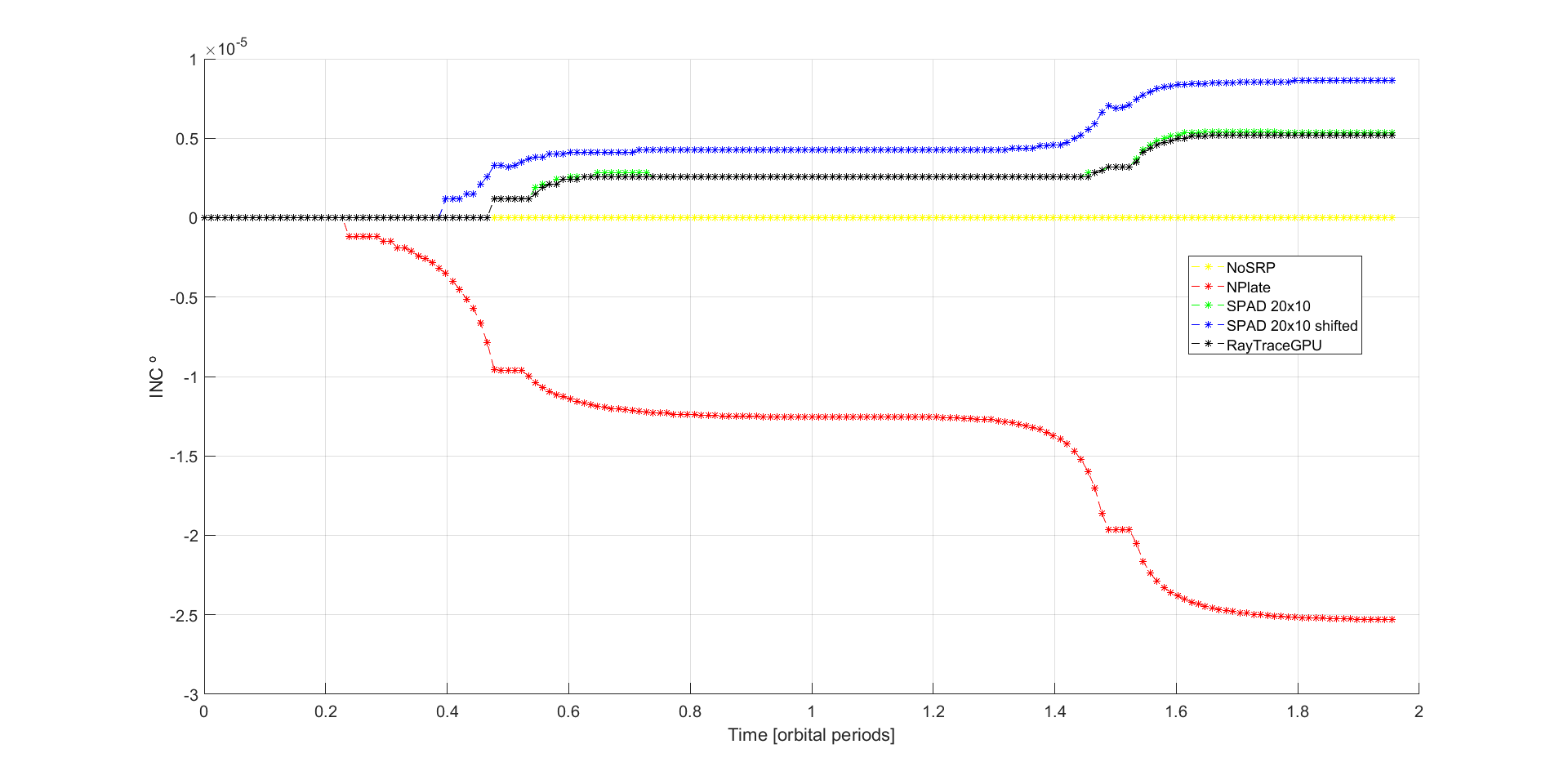}
\caption{Evolution of the orbital inclination for scenario $S_2$-E(S).}
\label{fig:closest_normal_inclination}
\end{figure}

As a representative example, Figures~\ref{fig:closest_normal_3d}-\ref{fig:closest_normal_inclination} show the results for scenario $S_2$-E(S), where the spacecraft changes its orientation at an azimuth angle of $65.82^\circ$. Figure~\ref{fig:closest_normal_z} shows the evolution of the Z-axis position at aphelion and perihelion.

The position difference between RayTrace GPU and SPAD 20x10 remains around $1\times10^{-8}$ AU (approximately 1.5 km), demonstrating that the SPAD interpolation accurately reproduces the high-fidelity solution. In contrast, SPAD 20x10S presents deviations of approximately $3\times10^{-7}$ AU (44.9 km), while N-Plate reaches differences between 224 km and 374 km.

A similar trend is observed for scenarios $S_2$ and $S_3$ across all orbital configurations. In both cases, E(S) consistently exhibits the largest deviations due to the stronger SRP environment experienced near perihelion.

For scenario $S_1$, the differences between the SRP models remain comparatively small, indicating that the spacecraft attitude changes occurring in this configuration have a limited influence on the accumulated SRP effects.

\paragraph{Continuous Attitude Variation ($S_4$ - E(S)).}

\begin{figure}[htbp]
\centering
\includegraphics[width=\textwidth]{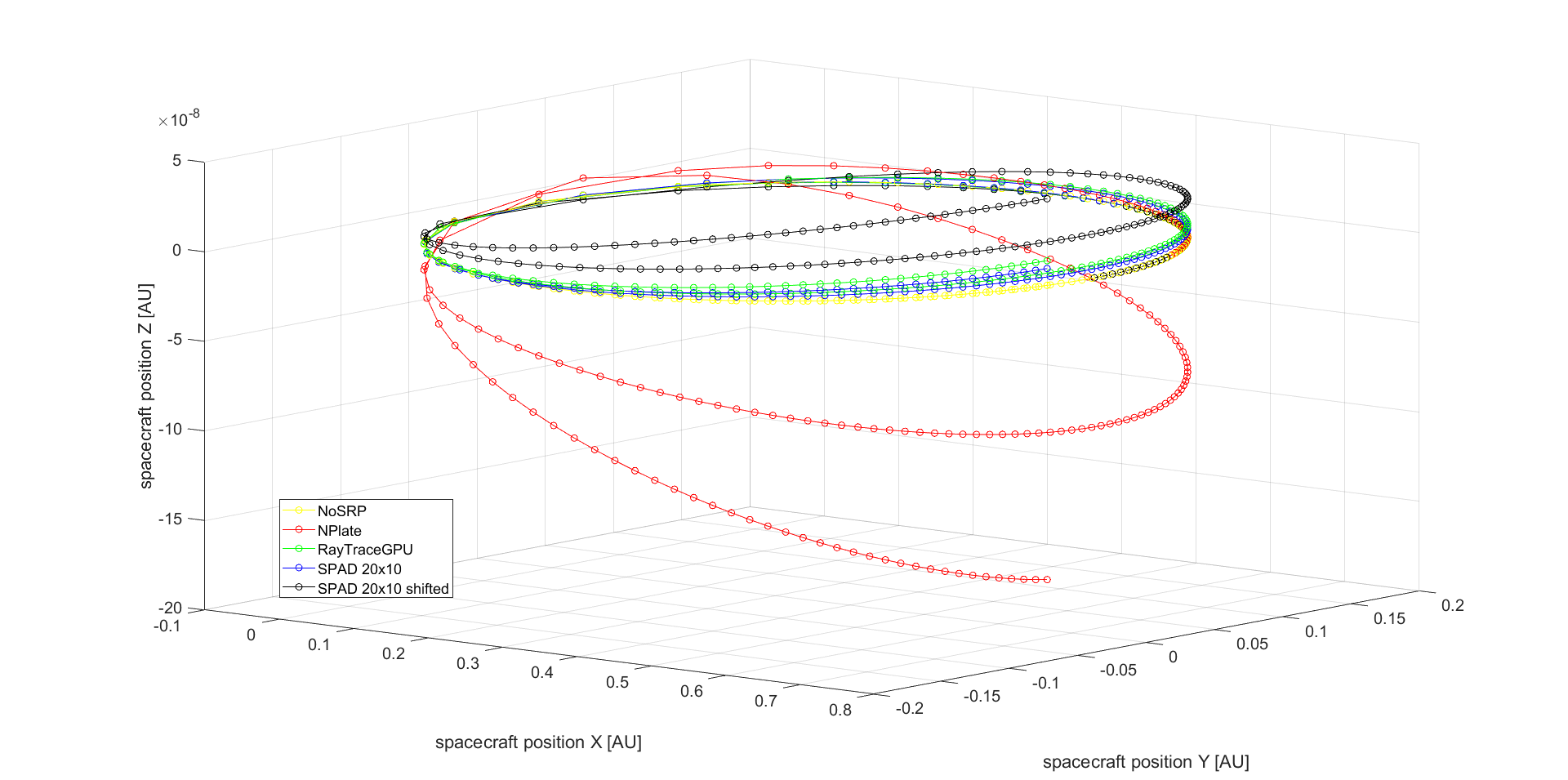}
\caption{Three-dimensional trajectory for the continuous-attitude scenario.}
\label{fig:continuous_3d}
\end{figure}

\begin{figure}[htbp]
\centering
\includegraphics[width=\textwidth]{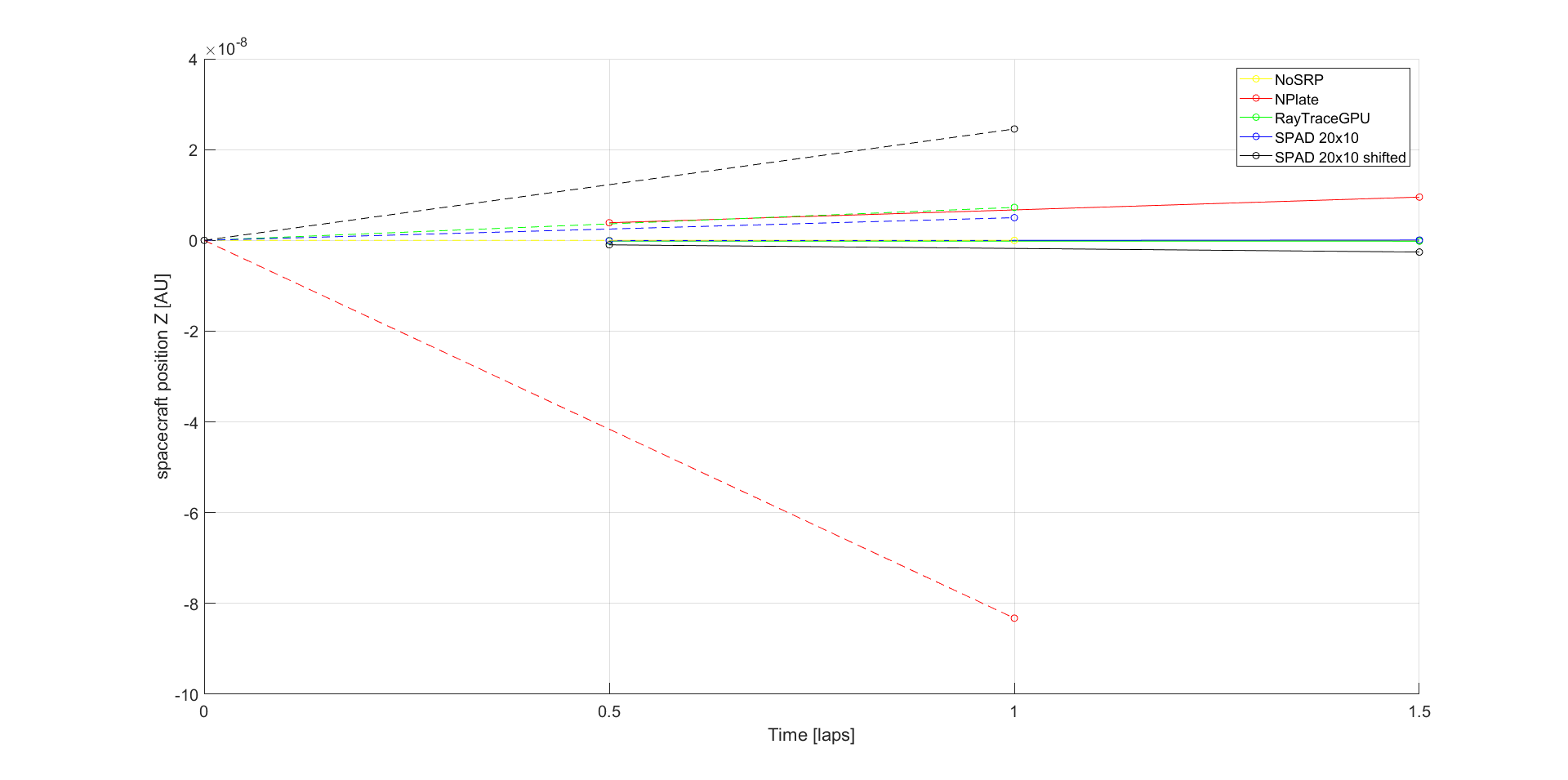}
\caption{Evolution of the spacecraft $z$-coordinate at aphelion and perihelion for the continuous-attitude scenario.}
\label{fig:continuous_z}
\end{figure}

\begin{figure}[htbp]
\centering
\includegraphics[width=\textwidth]{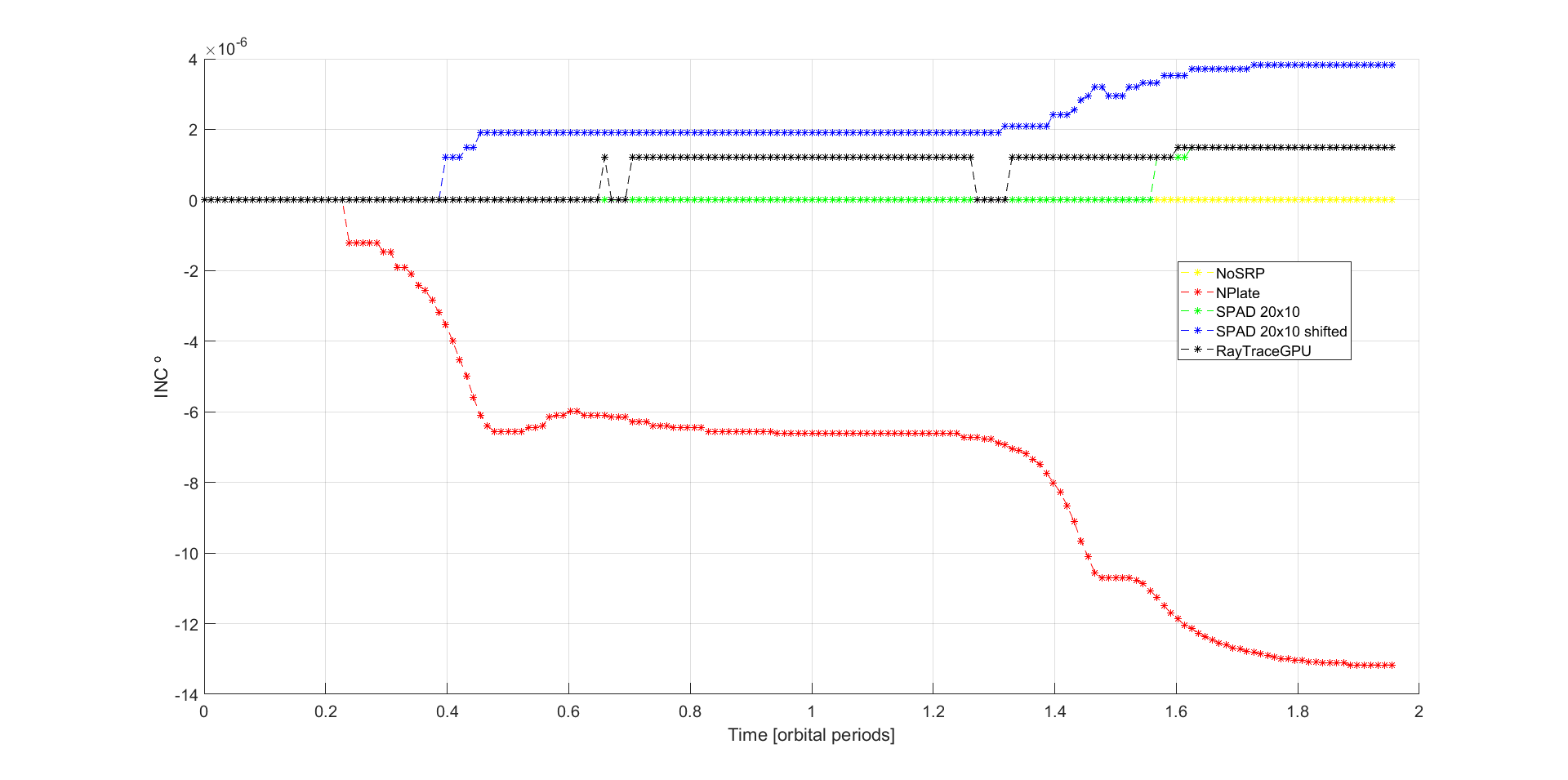}
\caption{Evolution of the orbital inclination for the continuous-attitude scenario.}
\label{fig:continuous_inclination}
\end{figure}

The continuous-attitude scenario is summarized in Figures~\ref{fig:continuous_3d}-\ref{fig:continuous_inclination}.

In this configuration, the spacecraft attitude varies gradually rather than switching between discrete orientations. As a result, the spacecraft spends longer periods in intermediate attitudes where secondary reflections contribute more significantly to the SRP force.

As shown in Figure~\ref{fig:continuous_z}, the deviation between RayTrace GPU and SPAD 20x10 reaches approximately $5.25\times10^{-9}$ AU (0.78 km). Although the error remains small, it is larger than in several of the discrete-attitude cases, suggesting that continuous attitude transitions are more sensitive to interpolation errors.

\paragraph{Computational Performance.}

Table~\ref{table:time_elliptic_orbits} summarizes the execution times for the family of elliptical orbit simulations.

\begin{table}[htbp]
\centering
\caption{Computational times of the elliptic orbits.}
\label{table:time_elliptic_orbits}
\small
\begin{tabular}{llccc}
\toprule
\multicolumn{2}{l}{\textbf{Computational Time}} & \textbf{$S_3$} & \textbf{$S_2$} & \textbf{$S_1$} \\
\midrule

\multirow{5}{*}{E(S)} 
& No SRP       & 58.054 ms  & 76.287 ms  & 65.478 ms  \\
& N-Plate       & 70.392 ms  & 70.146 ms  & 89.899 ms  \\
& SPAD 20x10   & 98.744 ms  & 120.360 ms & 93.123 ms  \\
& SPAD 20x10S  & 93.319 ms  & 98.481 ms  & 97.208 ms  \\
& RayTrace GPU & 56.93 min  & 1.05 h     & 1.09 h     \\
\midrule

\multirow{5}{*}{E(M)} 
& No SRP       & 71.345 ms  & 69.815 ms  & 68.436 ms  \\
& N-Plate       & 77.446 ms  & 90.252 ms  & 98.209 ms  \\
& SPAD 20x10   & 111.747 ms & 107.535 ms & 98.596 ms  \\
& SPAD 20x10S  & 102.083 ms & 114.376 ms & 101.416 ms \\
& RayTrace GPU & 1.09 h     & 1.25 h     & 1.34 h     \\
\midrule

\multirow{5}{*}{E(L)} 
& No SRP       & 108.896 ms & 82.178 ms  & 89.260 ms  \\
& N-Plate       & 103.113 ms & 101.090 ms & 100.415 ms \\
& SPAD 20x10   & 135.244 ms & 128.206 ms & 140.680 ms \\
& SPAD 20x10S  & 130.796 ms & 131.230 ms & 142.050 ms \\
& RayTrace GPU & 1.54 h     & 1.81 h     & 1.97 h     \\
\bottomrule
\end{tabular}
\end{table}

The results show that all interpolation-based approaches complete in less than one second, whereas direct RayTrace GPU propagation requires between approximately one and two hours depending on the orbital scenario.

This difference corresponds to several orders of magnitude in computational cost. Considering that the trajectory deviations between RayTrace GPU and SPAD 20x10 remain at kilometre-level or below for most scenarios, SPAD interpolation provides a highly favourable accuracy-to-performance ratio.

In particular, when the spacecraft attitude evolution is compatible with the SPAD discretization, the SPAD 20x10 model is able to reproduce the RayTrace GPU solution with very small trajectory deviations while reducing the computational time from hours to milliseconds.

\subsection{Discussion}

$RQ1$ investigated whether SPAD-based interpolation can replace direct RayTrace GPU evaluations during orbit propagation while preserving sufficient accuracy and providing a significant reduction in computational cost.

From the accuracy perspective, the results demonstrate that SPAD-based propagation reproduces the trajectories generated by direct RayTrace GPU evaluations with small position and orbital-element deviations for all scenarios considered. Among the interpolated approaches, SPAD 20x10 consistently provides the closest approximation to the RayTrace GPU reference. This behavior is explained by the characteristics of the analyzed scenarios, in which the spacecraft attitude varies primarily in azimuth while the elevation remains close to $0^\circ$. Since the SPAD 20x10 dataset contains samples at $0^\circ$ elevation, the interpolation is performed almost entirely along the azimuth direction, resulting in errors that are nearly equivalent to those of a one-dimensional interpolation. In contrast, the shifted SPAD 20x10S dataset does not contain samples at $0^\circ$ elevation, forcing the SRP acceleration to be estimated through bilinear interpolation between neighbouring samples. Consequently, this dataset provides a more representative assessment of the behavior of the bilinear interpolator itself. Although larger interpolation errors are observed, both SPAD-based approaches remain substantially more accurate than the N-Plate model.

Regarding computational performance, the improvement provided by SPAD interpolation is considerable. While direct RayTrace GPU propagation requires between one and two hours depending on the scenario, SPAD-based propagation completes in only a few milliseconds. This reduction of several orders of magnitude makes interpolation-based propagation more suitable for long-term trajectory simulations.

The results obtained with the shifted SPAD 20x10S dataset provide a more representative assessment of the bilinear interpolation strategy. Although the interpolation error increases when the reference elevation samples are no longer available, the resulting trajectories remain considerably closer to the direct RayTrace GPU solution than those obtained with the N-Plate model. More general attitude profiles involving simultaneous azimuth and elevation variations are therefore expected to increase the interpolation error and may require denser SPAD discretizations or, in the most demanding cases, direct high-fidelity SRP evaluations.

Overall, the results indicate that precomputed SPAD datasets provide an excellent compromise between accuracy and computational efficiency for the orbit propagation scenarios analyzed in this work. Direct RayTrace GPU evaluations remain essential for generating the reference SPAD datasets and for applications requiring the highest possible SRP fidelity, but their computational cost currently limits their practical use for routine orbit propagation.

\section{Acceleration of Direct High-Fidelity SRP Evaluation Using a Vulkan-Based Strategy}

This section evaluates the proposed Vulkan-based implementation using the experimental methodology described in the previous section. The results compare the Vulkan and OpenGL implementations in terms of numerical consistency and computational performance. First, the numerical equivalence of both implementations is verified to ensure that the proposed GPU backend preserves the original SRP formulation. Subsequently, the computational cost of individual SRP evaluations and complete orbit propagation simulations is analyzed to quantify the performance improvements provided by the Vulkan implementation.

\subsection{Numerical Validation of the Vulkan Implementation}

Before analyzing its computational performance, the numerical consistency of the proposed Vulkan implementation is verified against the original OpenGL-based ray-tracing implementation. Since both implementations solve the same SRP formulation, any differences should be limited to floating-point arithmetic and implementation details.

The validation is performed using the Parker Solar Probe spacecraft model described in Table~\ref{table:datasetsVulkanPM}. The SRP acceleration computed by the Vulkan implementation is compared with the reference OpenGL implementation for the 1,296 illumination directions defined in the experimental methodology.

The maximum relative discrepancy remains below $5\times10^{-4}$, confirming that both implementations produce numerically consistent results for practical SRP computations. Consequently, the accuracy metrics ($A_1$ and $A_2$) are not further analyzed, as both implementations solve the same SRP formulation and produce indistinguishable trajectories under the evaluated conditions. The remaining evaluation therefore focuses exclusively on the performance metrics ($P_1$ and $P_2$), allowing any observed differences to be attributed solely to the GPU implementation.

\subsection{Results}

Metric $P_1$ evaluates the computational cost of individual SRP force evaluations. Three implementations are compared: the original OpenGL-based RayTrace GPU implementation, the Vulkan implementation with CPU-based reduction, and the Vulkan implementation with GPU-based reduction. In the latter, the SRP contributions are accumulated directly on the GPU before transferring only the final result to the CPU.

Table~\ref{tab:P1} summarizes both the total computation time and the mean execution time per SRP evaluation.

\begin{table}[htbp]
\centering
\caption{Summary table of P1 results.}
\label{tab:P1}
\small
\begin{tabular}{lcc}
\toprule
 & \textbf{Total time} & \textbf{Mean time} \\
\midrule
PSP          & 359.499s & 0.277s \\
PSP Vulkan (sum GPU)   & 38.387s  & 0.028s \\
PSP Vulkan (sum CPU)  & 68.224s  & 0.049s \\
BoxWing          & 15.765s  & 0.012s \\
BoxWing Vulkan (sum GPU)  & 3.581s  & 0.003s \\
BoxWing Vulkan (sum CPU)  & 24.335s  & 0.017s \\
\bottomrule
\end{tabular}
\end{table}

The observed performance can be interpreted by considering the total execution time as the sum of a computation-dependent component and several overhead terms:
\begin{equation}
T_{\mathrm{total}} = T_{\mathrm{compute}} +T_{\mathrm{transfer}} +T_{\mathrm{reduction}},\label{eq:performance_model}\end{equation}

where $T_{\mathrm{compute}}$ corresponds to the ray tracing computation, $T_{\mathrm{transfer}}$ represents the time required to transfer intermediate results between the GPU and the CPU, and $T_{\mathrm{reduction}}$ denotes the accumulation of the SRP contributions.

In the GPU reduction approach, only the final accumulated force is transferred back to the CPU, making $T_{\mathrm{transfer}}$ negligible. Conversely, when the reduction is performed on the CPU, the complete buffer of intermediate contributions must be transferred from GPU memory before the accumulation can be carried out, increasing both the communication overhead and the total execution time.

The Vulkan implementation with GPU reduction provides the best performance for both spacecraft models. For the PSP model, the execution time decreases from 359.499\,s to 38.387\,s, corresponding to a speed-up of approximately 9.4. For the BoxWing model, it decreases from 15.765\,s to 3.581\,s, achieving a speed-up of approximately 4.4. These improvements result from performing both the ray tracing and the reduction entirely on the GPU, thereby minimizing communication overhead.

When the reduction is performed on the CPU, the performance becomes strongly dependent on the complexity of the spacecraft geometry. For the PSP model, the Vulkan implementation still reduces the total execution time to 68.224\,s, representing a speed-up of approximately 5.3 compared with the original OpenGL implementation. In contrast, for the BoxWing model, the CPU-based reduction increases the execution time from 15.765\,s to 24.335\,s, making the Vulkan implementation slower than the original OpenGL approach.

This behavior can be explained by the relative contribution of the different terms in Equation~\ref{eq:performance_model}. While $T_{\mathrm{transfer}}$ and $T_{\mathrm{reduction}}$ remain approximately constant because the output buffer size is identical in both experiments, $T_{\mathrm{compute}}$ increases with the complexity of the spacecraft geometry.
For the BoxWing spacecraft, which consists of only 212 triangles, the ray tracing workload is relatively small. Consequently, the communication and reduction overhead dominate the total execution time, offsetting the computational advantages of the Vulkan implementation. Under these conditions, the mature optimizations available in modern OpenGL graphics drivers, together with the relatively small computational workload, allow the original graphics-pipeline implementation to outperform the Vulkan compute-based approach in this case, as the fixed synchronization and data-transfer costs are not sufficiently amortized. This behavior is consistent with the design philosophy of Vulkan, which exposes resource management, memory allocation, and synchronization explicitly to the application instead of relying on implicit driver optimizations, as is the case in traditional APIs such as OpenGL \cite{VulkanProgrammingGuide}.

In contrast, the PSP spacecraft contains approximately 39,000 triangles, making the ray tracing computation substantially more expensive. As the computational workload increases, the fixed communication overhead becomes progressively less significant relative to the total execution time. Consequently, the improved computational efficiency of the Vulkan implementation outweighs the additional transfer cost, resulting in a significant reduction in execution time even when the reduction is performed on the CPU.

The impact of the Vulkan implementation on complete orbit propagation is evaluated through metric $P_2$, which measures the total computation time required to simulate 20 orbital revolutions. Based on the results for metric $P_1$, only the GPU reduction strategy is considered for the orbit propagation experiments.

The results are summarized in Table~\ref{tab:P2_Vulkan}.

\begin{table}[htbp]
\centering
\caption{Summary table of P2 results.}
\label{tab:P2_Vulkan}
\small
\begin{tabular}{lc}
\toprule
 & \textbf{Total time (min)} \\
\midrule
PSP             & 318.05 \\
PSP Vulkan      & 20.88 \\
BoxWing         & 12.92 \\
BoxWing Vulkan  & 3.13 \\
\bottomrule
\end{tabular}
\end{table}

The results demonstrate that the performance improvements observed for individual SRP evaluations are maintained when integrated into a complete orbit propagation framework. For the PSP spacecraft, the total propagation time is reduced from 5h 18min to 20min 53s, corresponding to a speed-up of approximately 15.2. Similarly, for the BoxWing spacecraft, the execution time decreases from 12min 55s to 3min 8s, achieving a speed-up of approximately 4.1.

The larger speed-up obtained for the PSP model is consistent with the observations from metric~$P_1$. Since orbit propagation requires thousands of SRP evaluations, the computational savings achieved for each evaluation accumulate throughout the propagation process, particularly for complex spacecraft geometries where SRP computation dominates the execution time.

Overall, the results show that the proposed Vulkan implementation substantially improves the scalability of the SRP simulation framework while preserving the numerical fidelity of the original RayTrace GPU model. This enables the practical use of direct high-fidelity SRP evaluation in long-duration orbit propagation scenarios where the computational cost of the original implementation would otherwise be prohibitive.

\subsection{Discussion}

From the numerical accuracy perspective, the proposed Vulkan implementation preserves the numerical fidelity of the original OpenGL-based RayTrace GPU model. The relative differences in the computed SRP accelerations remain below $5\times10^{-4}$, confirming that the optimization modifies only the computational implementation without affecting the underlying physical model.

From the computational performance perspective, the proposed implementation substantially reduces the cost of direct high-fidelity SRP evaluation. The largest improvements are obtained when the reduction stage is performed entirely on the GPU, eliminating the transfer of intermediate buffers to the CPU. Under this configuration, the execution time of individual SRP evaluations is reduced by up to a factor of 9.4 for the Parker Solar Probe model. These improvements are directly reflected in complete orbit propagation simulations, where the total propagation time is reduced by up to a factor of 15.2.

The performance gains depend on the computational complexity of the spacecraft model. For simple geometries, such as the BoxWing spacecraft, the ray-tracing workload is relatively small and communication overhead represents a significant fraction of the total execution time. In contrast, for complex spacecraft such as the Parker Solar Probe, the ray-tracing computation dominates the execution time, allowing the Vulkan implementation to fully exploit the available GPU parallelism and achieve substantially larger speed-ups.

Overall, these results answer RQ2 by demonstrating that direct high-fidelity SRP evaluation can be significantly accelerated through the proposed Vulkan implementation without compromising numerical accuracy. The proposed approach therefore makes high-fidelity SRP computation considerably more practical for long-duration orbit propagation while preserving the physical consistency of the original RayTrace GPU model.

\section{Increasing Realism in SRP Simulation by Using Movable Solar Panels}

This section evaluates the impact of incorporating movable solar panels into the SRP model using the experimental methodology defined in the previous section. The results assess both the improvement in orbit propagation accuracy and the associated computational overhead with respect to the fixed-panel configuration.

Before presenting the quantitative results, the spacecraft configurations with movable panels used in the experiments are briefly described. For the Parker Solar Probe (PSP), both solar arrays rotate about their attachment axis, allowing limited forward and backward motion within an angular range of approximately $\pm30^\circ$. This configuration represents the partial retraction of the arrays during operation. For the BoxWing spacecraft, each solar array is modeled with two rotational degrees of freedom. The first allows the panel to rotate about its support axis within approximately $\pm90^\circ$, enabling Sun-pointing. The second models the hinge rotation of the panel, allowing an additional folding motion within approximately $\pm45^\circ$. Together, these rotations provide a simplified representation of a two-axis movable solar-array mechanism.

\begin{figure}[htbp] \centering \begin{subfigure}[t]{0.35\textwidth} \centering \includegraphics[width=\textwidth]{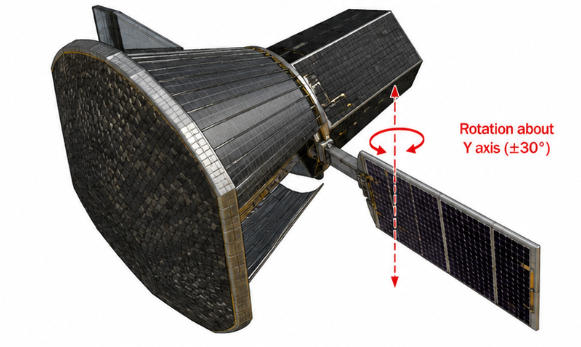} \\[2mm] \caption{PSP} \end{subfigure} \hfill \begin{subfigure}[t]{0.64\textwidth} \centering \includegraphics[width=\textwidth]{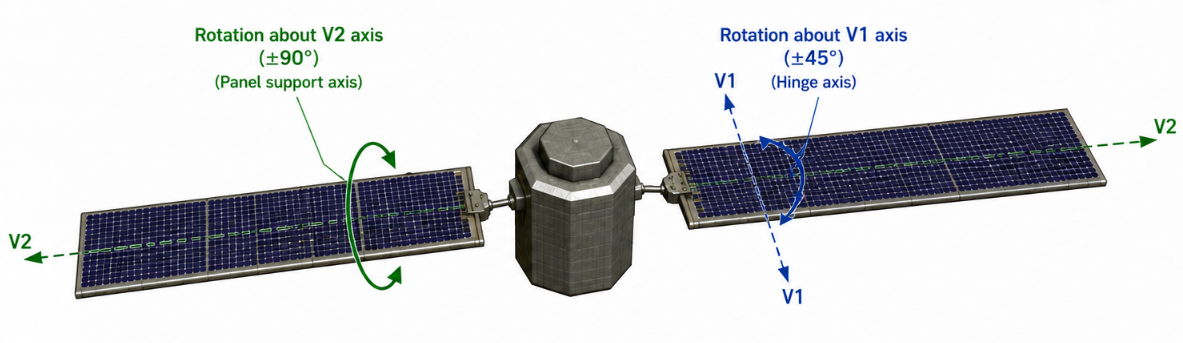} \\[2mm] \caption{BoxWing} \end{subfigure} \caption{Movable solar panel configurations considered for the Parker Solar Probe and BoxWing spacecraft models.} \label{fig:mobile_panels} \end{figure}

Figure~\ref{fig:mobile_panels} illustrates the movable configurations considered for both spacecraft models.

\subsection{Qualitative Validation of the Movable Solar Panel Strategy}

Although no ground-truth dataset is available to validate the movable solar-panel model, the correctness of the implementation can be assessed qualitatively. The proposed approach only modifies the orientation of the movable panels by updating the positions of their vertices and the corresponding surface normals according to the prescribed rotation. The ray tracing algorithm, the SRP computation, and the GPU implementation remain identical to those validated for the static spacecraft model.

Since only the panel vertex positions and normals are modified while the ray tracing algorithm remains unchanged, verifying that the articulated geometry is correctly generated is sufficient to validate the implementation of the geometric update.

\begin{figure}[htbp]
\centering

\begin{subfigure}[b]{0.48\textwidth}
    \centering
    \includegraphics[width=\textwidth]{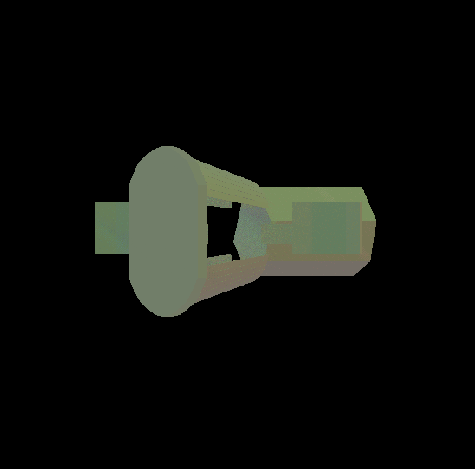}
    \caption{PSP spacecraft.}
    \label{fig:PSP_PM_validation}
\end{subfigure}
\hfill
\begin{subfigure}[b]{0.427\textwidth}
    \centering
    \includegraphics[width=\textwidth]{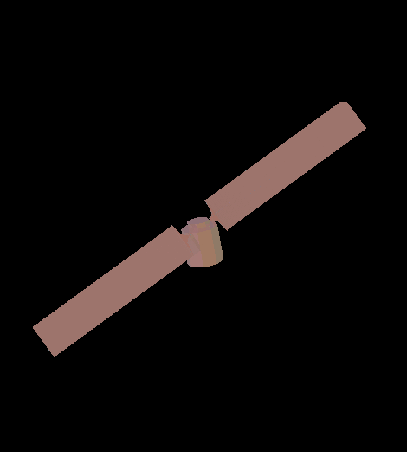}
    \caption{BoxWing spacecraft.}
    \label{fig:BoxWing_PM_validation}
\end{subfigure}

\caption{Qualitative validation of the movable solar-panel implementation. The rendered GPU output confirms that the solar panels are correctly rotated while the remainder of the spacecraft geometry is preserved.}
\label{fig:PM_validation}

\end{figure}

Consequently, visual inspection of the rendered spacecraft geometry from the GPU provides an effective means of verifying the implementation. Figure~\ref{fig:PM_validation} shows representative rendering results for the PSP and BoxWing spacecraft after applying panel rotations. The rendered images confirm that the movable panels rotate as expected while preserving the remainder of the spacecraft geometry, providing confidence that the updated geometry supplied to the SRP computation is correct.

\subsection{Results}

This section presents the experimental results obtained with the movable solar panel model. The fixed- and movable-panel configurations are compared in terms of SRP force accuracy, orbit propagation, and computational performance. 

The results for metric $A_1$ are summarized in Table~\ref{tab:A1}. They quantify the error introduced by neglecting movable solar panels in the computation of the SRP forces.

\begin{table}[htbp]
\centering
\caption{Summary table of A1 results.}
\label{tab:A1}
\small
\begin{tabular}{lccc}
\toprule
 &  \textbf{MSE} & \textbf{max}  & \textbf{MSER}\\
\midrule
PSP & 3.999e-11 & 1.289e-10 & 7.489e-02\\
BoxWing & 2.510e-11 & 4.648e-11 & 2.750\\
\bottomrule
\end{tabular}
\end{table}

For the PSP model, the mean relative error (MSER) when movable panels are neglected is 7.489e-02 (around 7.5\%). This indicates that neglecting panel motion introduces a moderate error in the computed SRP acceleration. This behavior is consistent with the compact geometry and relatively low area-to-mass ratio of the PSP spacecraft.

In contrast, the BoxWing satellite exhibits a much higher sensitivity, with a mean relative error of 2.750. This result indicates that neglecting panel motion introduces a substantially larger error for box-wing spacecraft, where the large movable solar arrays have a much greater influence on the resulting SRP acceleration.

While metric~$A_1$ quantifies the instantaneous SRP differences, metric~$A_2$ evaluates how these differences accumulate during long-term orbit propagation.

The reference elliptical orbit was propagated for 20 orbital revolutions for both spacecraft models, considering configurations with and without movable solar panels to evaluate metric~$A_2$.

\begin{figure}[ht!]

\makebox[\textwidth][c]{%
\begin{minipage}{1.05\textwidth}
    \centering
\includegraphics[width=\textwidth]{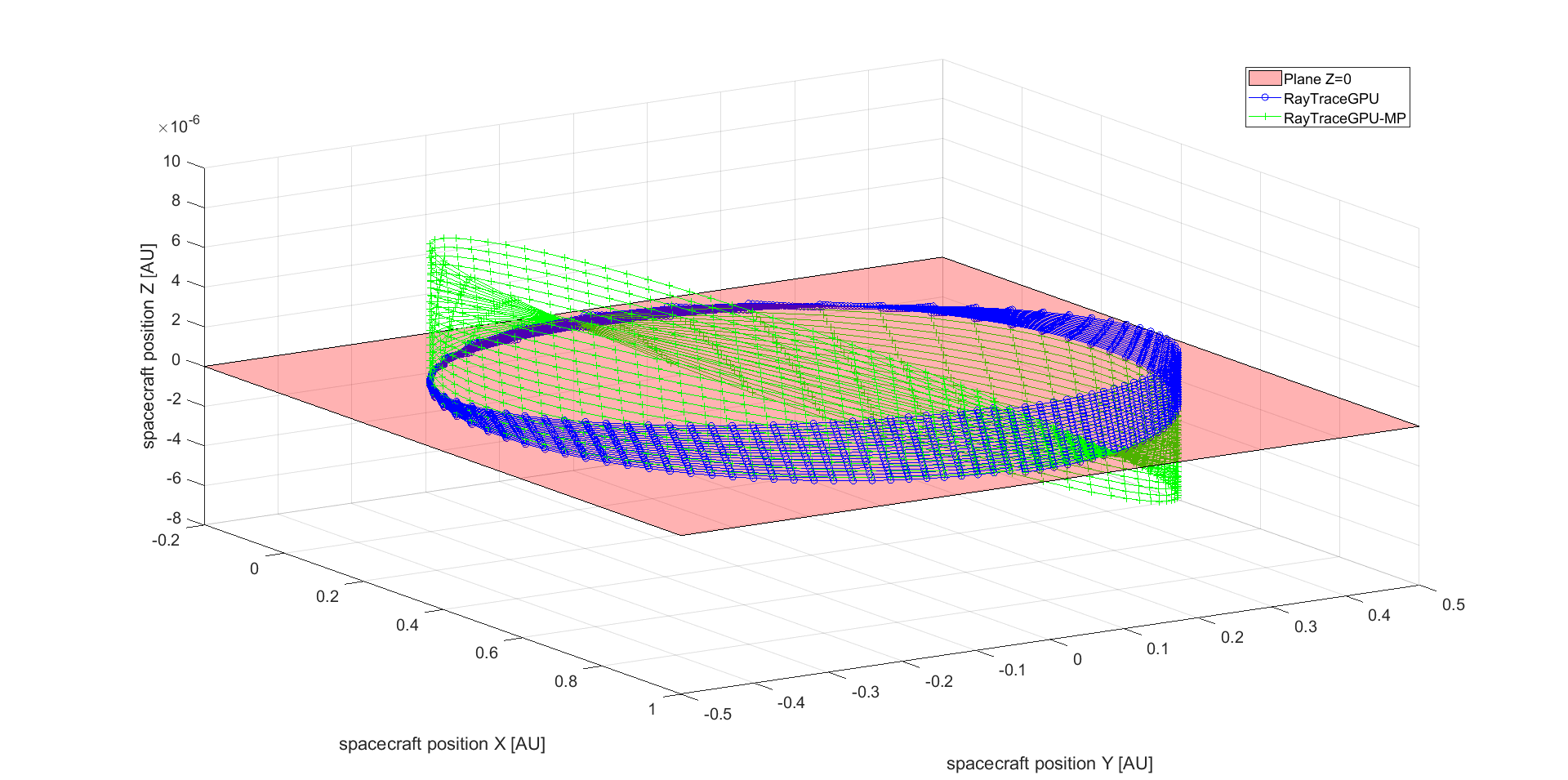}
\end{minipage}
}
\caption{3D view of the PSP trajectory using fixed and movable solar panels. The label RayTraceGPU corresponds to the fixed-panel configuration, whereas RayTraceGPUMP denotes the movable-panel configuration.}
\label{fig:PSP-XYZ}
\end{figure}

\begin{figure}[ht!]

\makebox[\textwidth][c]{%
\begin{minipage}{1.2\textwidth}
    \centering
\includegraphics[width=\textwidth]{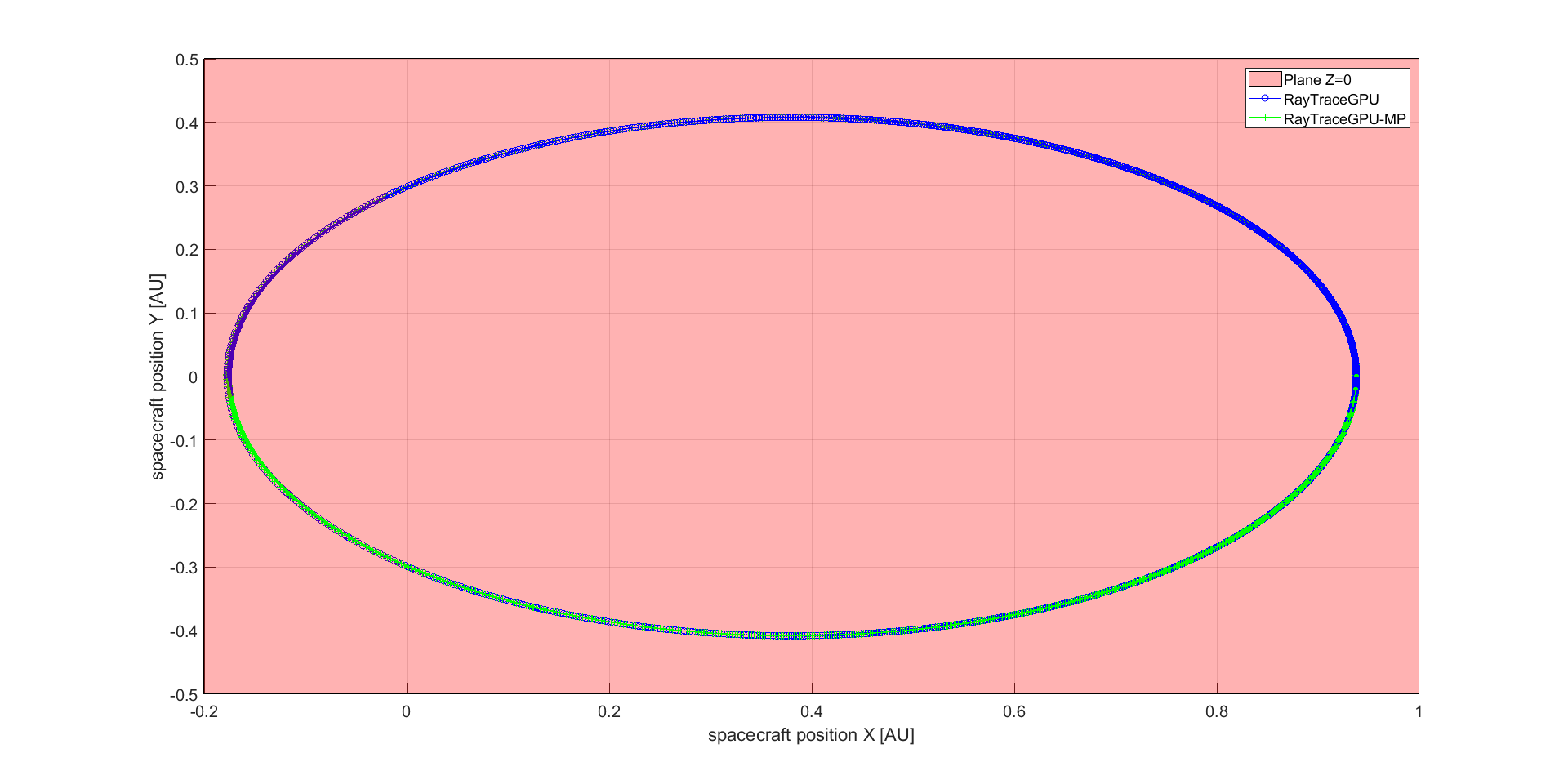}
\end{minipage}
}
\caption{XY projection of the PSP trajectory. The label RayTraceGPU corresponds to the fixed-panel configuration, whereas RayTraceGPUMP denotes the movable-panel configuration.}
\label{fig:PSP-XY}
\end{figure}

\begin{figure}[ht!]

\makebox[\textwidth][c]{%
\begin{minipage}{1.2\textwidth}
    \centering
\includegraphics[width=\textwidth]{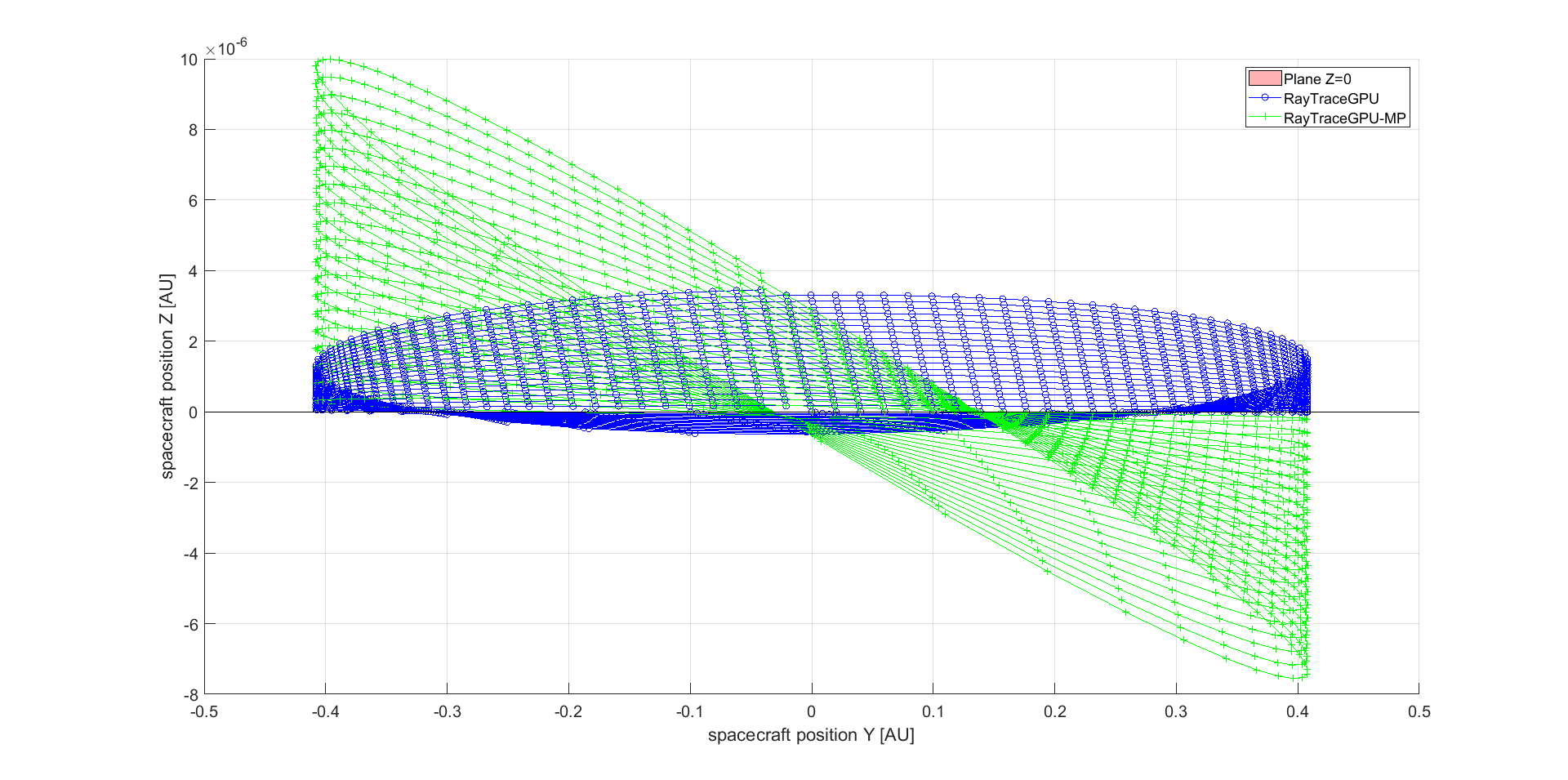}
\end{minipage}
}
\caption{YZ projection of the PSP trajectory. The label RayTraceGPU corresponds to the fixed-panel configuration, whereas RayTraceGPUMP denotes the movable-panel configuration.}
\label{fig:PSP-YZ}
\end{figure}

\begin{figure}[ht!]

\makebox[\textwidth][c]{%
\begin{minipage}{1.2\textwidth}
    \centering
\includegraphics[width=\textwidth]{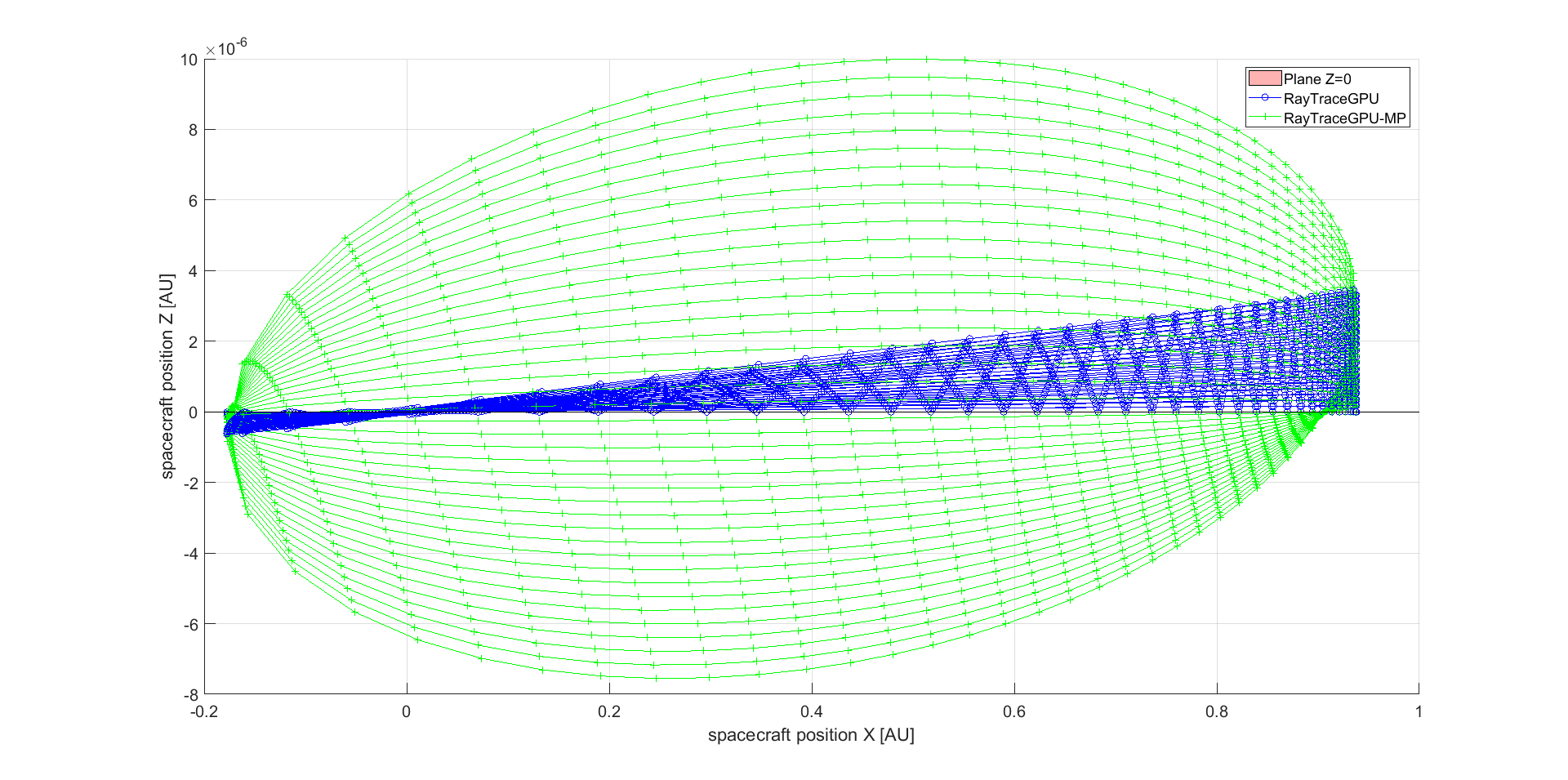}
\end{minipage}
}
\caption{XZ projection of the PSP trajectory. The label RayTraceGPU corresponds to the fixed-panel configuration, whereas RayTraceGPUMP denotes the movable-panel configuration.}
\label{fig:PSP-XZ}
\end{figure}

Figure~\ref{fig:PSP-XYZ} shows the resulting PSP trajectories, while Figures~\ref{fig:PSP-XY}, \ref{fig:PSP-YZ}, and \ref{fig:PSP-XZ} present their projections onto the principal planes. Although the instantaneous SRP differences quantified in metric~$A_1$ are moderate, their cumulative effect over multiple orbital revolutions produces a visible trajectory divergence.

\begin{figure}[ht!]

\makebox[\textwidth][c]{%
\begin{minipage}{1.1\textwidth}
    \centering
\includegraphics[width=\textwidth]{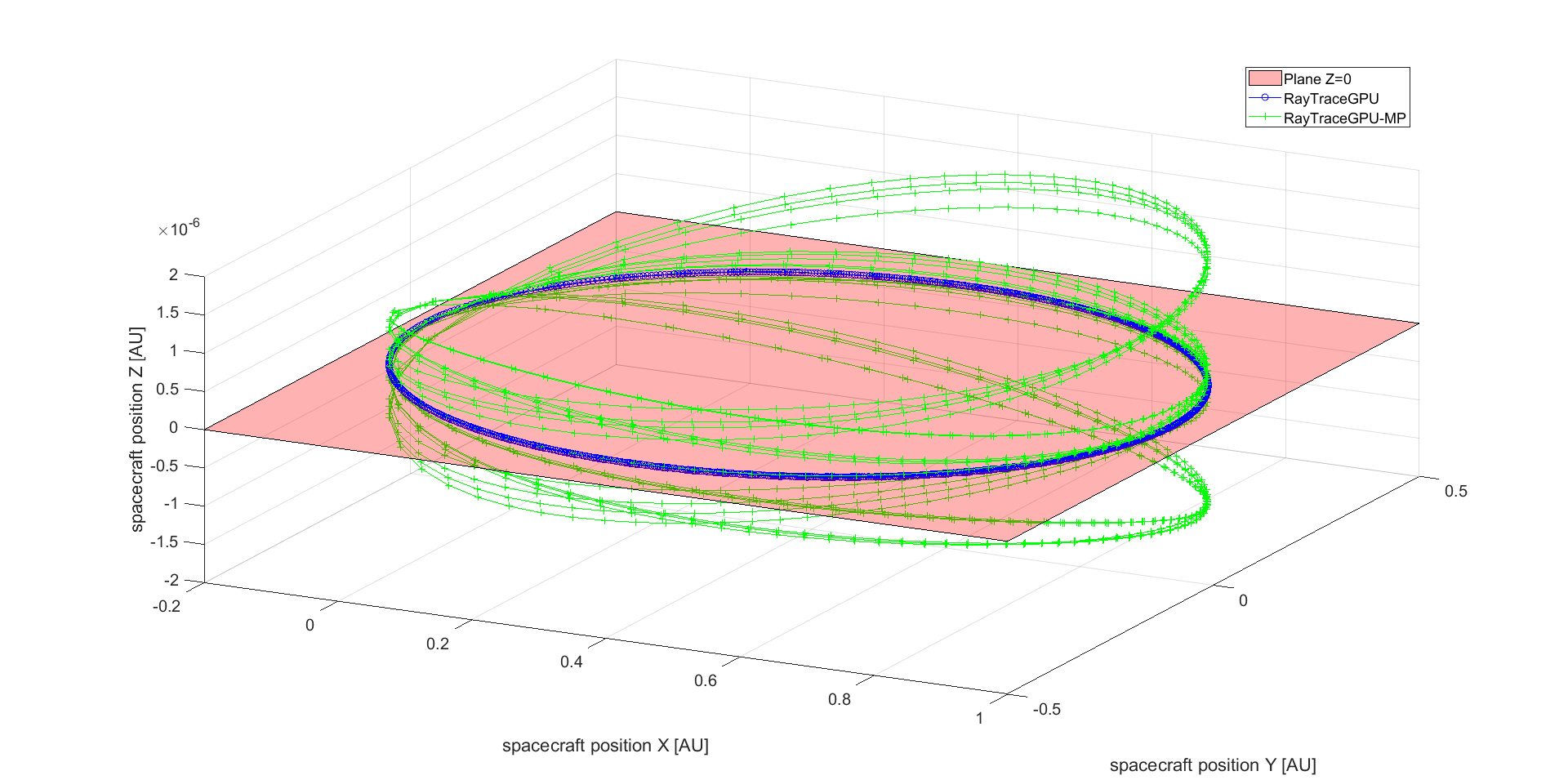}
\end{minipage}
}

\caption{3D view of the BoxWing trajectory using both static model and mobile panels. The label RayTraceGPU corresponds to the fixed-panel configuration, whereas RayTraceGPUMP denotes the movable-panel configuration.}
\label{fig:GPS-XYZ}
\end{figure}

\begin{figure}[ht!]

\label{fig:GPS-XY}
\makebox[\textwidth][c]{%
\begin{minipage}{1.2\textwidth}
    \centering
\includegraphics[width=\textwidth]{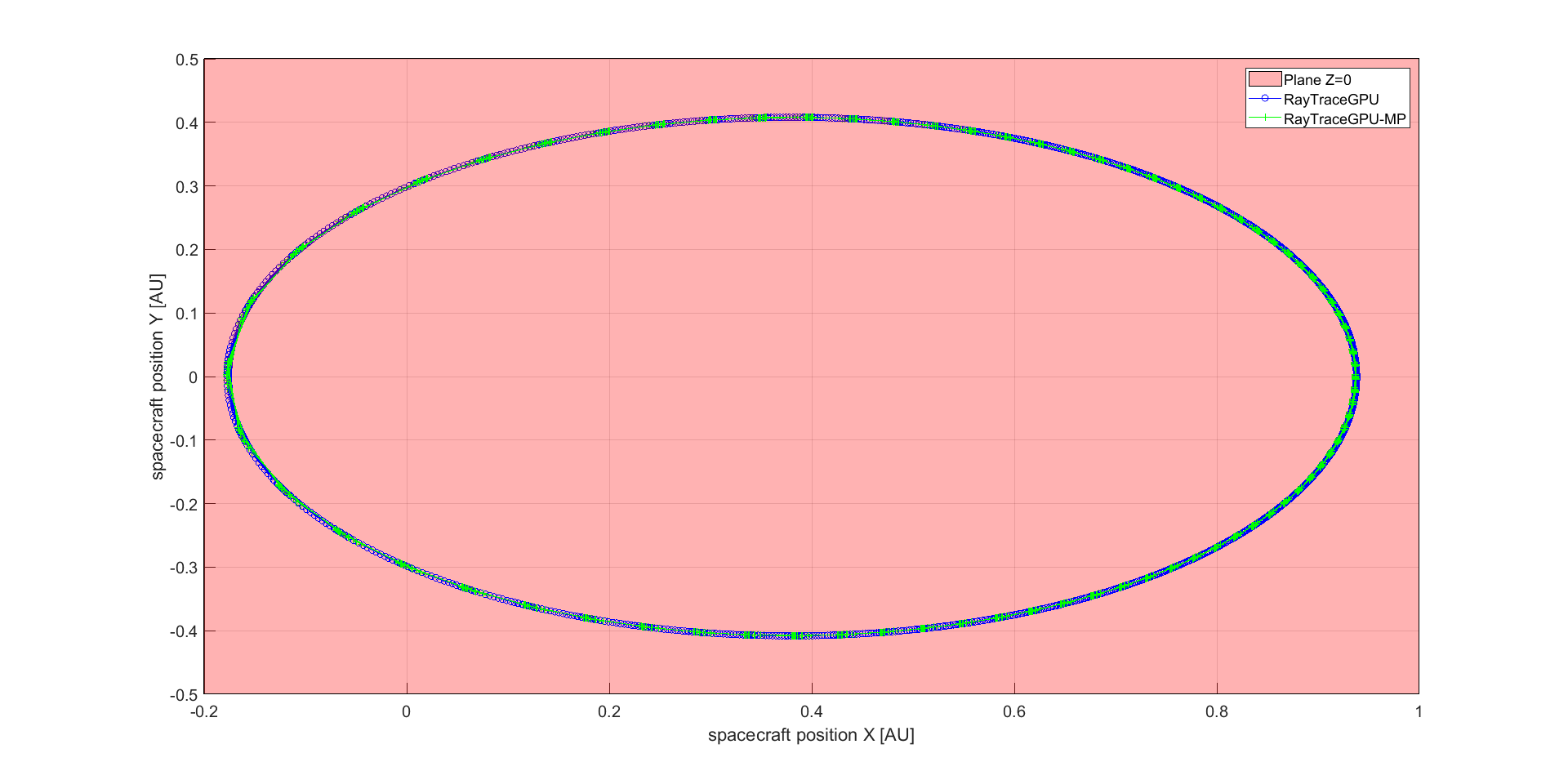}
\end{minipage}
}
\caption{XY projection of the BoxWing trajectory. The label RayTraceGPU corresponds to the fixed-panel configuration, whereas RayTraceGPUMP denotes the movable-panel configuration.}
\end{figure}

\begin{figure}[ht!]

\label{fig:GPS-YZ}
\makebox[\textwidth][c]{%
\begin{minipage}{1.2\textwidth}
    \centering
\includegraphics[width=\textwidth]{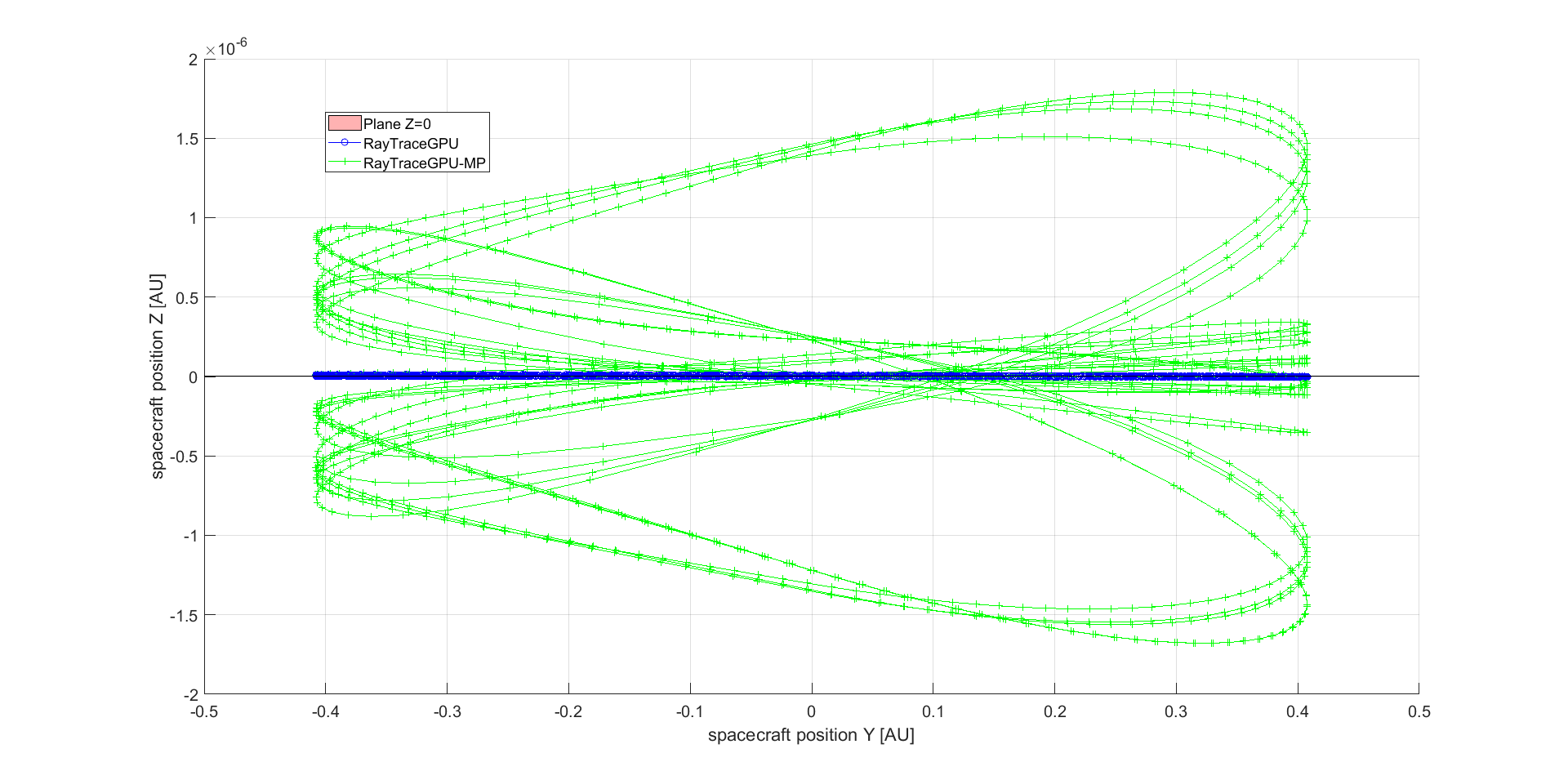}
\end{minipage}
}
\caption{YZ projection of the BoxWing trajectory. The label RayTraceGPU corresponds to the fixed-panel configuration, whereas RayTraceGPUMP denotes the movable-panel configuration.}
\end{figure}

\begin{figure}[ht!]

\makebox[\textwidth][c]{%
\begin{minipage}{1.2\textwidth}
    \centering
\includegraphics[width=\textwidth]{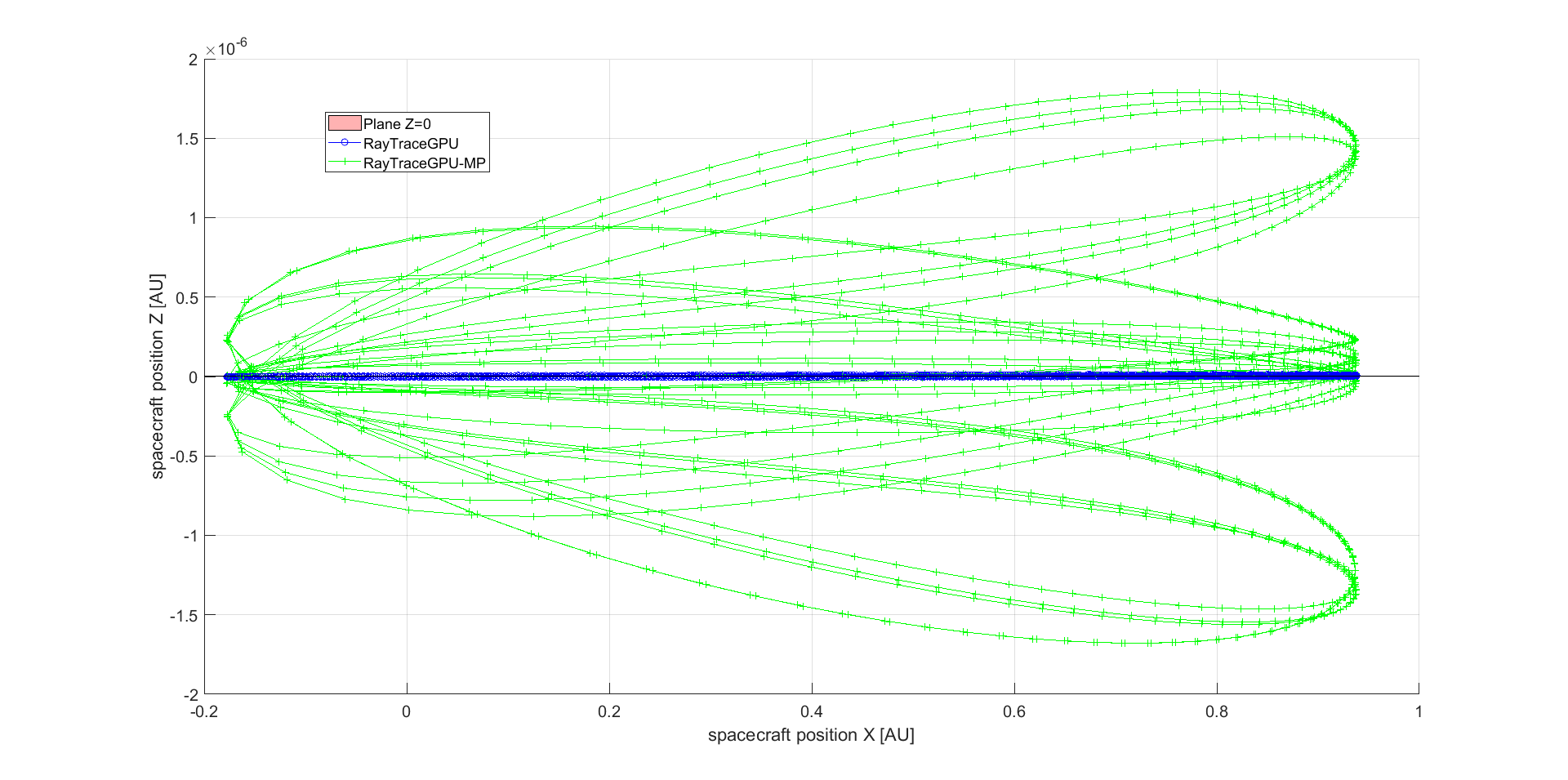}
\end{minipage}
}
\caption{XZ projection of the BoxWing trajectory. The label RayTraceGPU corresponds to the fixed-panel configuration, whereas RayTraceGPUMP denotes the movable-panel configuration.}
\label{fig:GPS-XZ}
\end{figure}

In contrast, Figures~\ref{fig:GPS-XYZ}-\ref{fig:GPS-XZ} illustrate the corresponding trajectories for the BoxWing spacecraft, where the larger SRP differences caused by panel articulation lead to a substantially greater orbital divergence.

Table~\ref{tab:A2} provides a quantitative summary of the trajectory deviations.

\begin{table}[htbp]
\centering
\caption{Summary table of A2 results [AU].}
\label{tab:A2}
\small
\begin{tabular}{lccc}
\toprule
 & \textbf{min} & \textbf{mean} & \textbf{max} \\
\midrule
PSP & 0.00 & 1.8443e-04 & 9.4028e-04 \\
BoxWing & 0.00 & 2.1175e-02 & 1.5668e-01 \\
\bottomrule
\end{tabular}
\end{table}

For the PSP spacecraft, the mean position difference between simulations with and without movable panels is approximately 27 590 km. Although the acceleration differences quantified in metric~$A_1$ are moderate, their cumulative effect over 20 revolutions results in a noticeable trajectory deviation.

For the BoxWing spacecraft, the impact is substantially larger. The mean trajectory deviation reaches approximately $3.167\times10^{6}$ km, reflecting the much stronger influence of panel orientation on the computed SRP acceleration. This behavior is consistent with the larger area-to-mass ratio and the box-wing geometry of the spacecraft.

These results demonstrate that even moderate differences in SRP acceleration can accumulate into significant orbital deviations over long propagation intervals. The effect becomes particularly pronounced for spacecraft with large movable solar arrays, highlighting the importance of accurately modeling panel orientation in high-fidelity SRP simulations.

The previous results quantify the accuracy improvements provided by movable solar panels. The following analysis evaluates the associated computational cost.

Metric $P_1$ measures the computational cost associated with SRP force evaluations for the 1,296 incident light directions in both the fixed- and movable-panel spacecraft configurations. In addition to the execution time, the total number of ray-triangle intersections is also analyzed, since it provides a direct measure of the ray-tracing workload.

\begin{table}[htbp]
\centering
\caption{Summary table of P1 results.}
\label{table:P1_PM}
\small
\begin{tabular}{lcc}
\toprule
 & \textbf{total time} & \textbf{mean time} \\
 \midrule
PSP     & 359.499s & 0.277s \\
PSP MP  & 353.506s & 0.258s \\
BoxWing     & 15.765s  & 0.012s \\
BoxWing MP  & 21.423s  & 0.015s \\
\bottomrule
\end{tabular}
\end{table}

Table~\ref{table:P1_PM} summarizes the total and mean computation times for the static and movable-panel configurations. To facilitate the interpretation of these results, Table~\ref{table:P1_PM_hits} also reports the total and mean number of ray-triangle intersections (hits), which are directly related to the amount of work performed during the ray tracing process.

\begin{table}[htbp]
\centering
\caption{Summary table of number of hits.}
\label{table:P1_PM_hits}
\small
\begin{tabular}{lcc}
\toprule & \textbf{Total hits} & \textbf{Mean hits} \\ 
\midrule
PSP     & 51 479 139 & 37 604 \\
PSP MP  & 50 503 997 & 36 891 \\
BoxWing     & 12 804 577 & 9 353 \\
BoxWing MP  & 25 811 091  & 18 854 \\
\bottomrule
\end{tabular}
\end{table}

For the PSP spacecraft, the movable-panel configuration produces a computation time that is comparable to the static configuration, with a slightly lower execution time despite the additional panel articulation. As shown in Table~\ref{table:P1_PM_hits}, the movable configuration also generates a slightly smaller number of ray-triangle intersections (approximately 2\% fewer hits on average). Consequently, the additional cost associated with updating the panel orientation has a negligible impact on the total execution time, while the reduced number of intersections slightly decreases the ray tracing workload.

A different behavior is observed for the BoxWing spacecraft. The movable-panel configuration almost doubles the average number of ray-triangle intersections, increasing from 9353 to 18854 hits per evaluation. This additional ray tracing workload is directly reflected in the computation time, which increases from 0.012\,s to 0.015\,s per SRP evaluation. Therefore, the increase in execution time is primarily caused by the larger number of ray-surface intersection tests rather than by the panel articulation itself.

Metric $P_2$ evaluates the computational cost of incorporating movable solar panels during complete orbit propagation by measuring the total execution time required to simulate 20 orbital revolutions.

\begin{table}[htbp]
\centering
\caption{Summary table of P2 results.}
\label{table:P2_PM}
\small
\begin{tabular}{lc}
\toprule
 & \textbf{Total time (min)} \\
\midrule
PSP        & 318.05 \\
PSP MP     & 343.42 \\
BoxWing    & 12.92 \\
BoxWing MP & 20.07 \\
\bottomrule
\end{tabular}
\end{table}

The total simulation time required to propagate the trajectory over 20 orbital revolutions is summarized in Table~\ref{table:P2_PM}. To facilitate the interpretation of the computational performance, Table~\ref{table:P2_PM_hits} also reports the total and mean number of ray-triangle intersections recorded during the propagation.

\begin{table}[htbp]
\centering
\caption{Summary table of number of hits during orbit propagation for both spacecraft models.}
\label{table:P2_PM_hits}
\small \begin{tabular}{lcc}
\toprule & \textbf{Total hits} & \textbf{Mean hits} \\
\midrule
 PSP     & 3\,707\,151\,720 & 44\,910\\
PSP MP  & 4\,059\,741\,820 & 49\,181 \\
BoxWing    &   891\,025\,343  & 10\,794 \\
BoxWing MP & 1\,441\,041\,286 & 17\,457 \\
\bottomrule
\end{tabular}
\end{table}

For both spacecraft models, the movable-panel configuration increases the total propagation time. For the PSP spacecraft, the execution time increases from 5h 18min (318.05 min) to 5h 43min (343.42 min), whereas for the BoxWing spacecraft it increases from 12.92 min to 20.07 min.

As shown in Table~\ref{table:P2_PM_hits}, this increase is consistent with the corresponding increase in the ray tracing workload. For the PSP spacecraft, the movable-panel configuration increases the average number of ray-triangle intersections from 44910 to 49181 per SRP evaluation (approximately 9.5\%). Similarly, for the BoxWing spacecraft, the average number of intersections increases from 10794 to 17457 (approximately 62\%). Since ray-triangle intersection tests constitute the dominant computational cost of the high-fidelity SRP model in these experiments, the larger number of intersections directly explains the increase in propagation time observed for both spacecraft.

This behavior differs from the results obtained for metric~$P_1$, where the PSP movable-panel configuration produced a slightly lower number of intersections than the static spacecraft. The difference arises because the two experiments evaluate different illumination scenarios. The experiment used for metric~$P_1$ considers a uniform distribution of incident light directions over the sphere, whereas orbit propagation only evaluates the sequence of Sun-spacecraft geometries encountered along the simulated orbit. Consequently, the average ray-triangle workload during propagation differs from that obtained using a uniform directional sampling.

The different behavior observed for the PSP and BoxWing spacecraft is likely related to their geometric characteristics. The BoxWing model is dominated by two large movable solar arrays, whose orientation has a pronounced influence on the number of ray-triangle intersections. Consequently, both the uniform directional sampling of the sphere for metric~$P_1$ and the illumination conditions encountered during orbit propagation produce the same overall trend. In contrast, the movable panels of the PSP spacecraft are considerably smaller and have a more limited range of motion relative to the spacecraft body. As a result, the average ray-triangle workload appears to be more sensitive to the particular illumination conditions encountered in each experiment. Although a more detailed sensitivity analysis would be required to quantify this effect, the observed results are consistent with the different geometric characteristics of both spacecraft.

Overall, the computational overhead associated with movable solar panels is primarily driven by the additional ray-triangle intersections generated by the changing spacecraft geometry. Despite this increase in ray-tracing workload, the execution-time overhead remains moderate relative to the improvement in physical realism provided by the high-fidelity SRP model.

\subsection{Discussion}

The incorporation of movable solar panels provides a more physically realistic representation of SRP during orbit propagation, particularly for spacecraft configurations with large movable solar arrays. The results show that neglecting panel mobility can introduce substantial errors in both instantaneous SRP force computation and long-term trajectory propagation.

Regarding SRP force accuracy, metric $A_1$ shows that the impact of neglecting panel motion depends strongly on the spacecraft geometry. For the PSP spacecraft, the relative error in the computed SRP force remains moderate, with an MSER of approximately 7.5\%. In contrast, the BoxWing spacecraft exhibits a significantly larger error, with an MSER of 2.75, due to the dominant contribution of its large movable solar arrays. These instantaneous differences accumulate over time, leading to noticeable trajectory deviations. After 20 orbital revolutions, the mean separation between fixed- and movable-panel trajectories reaches approximately 27,590 km for PSP and $3.17\times10^6$ km for the BoxWing spacecraft.

The inclusion of movable panels also introduces an additional computational cost. The performance results show that this overhead is mainly associated with the increased number of ray-triangle intersection tests caused by the changing spacecraft geometry rather than by the panel orientation updates themselves. For complete orbit propagation, the computational increase remains moderate for PSP, with the execution time increasing by approximately 8\%, whereas the BoxWing configuration presents a larger increase of approximately 55\% due to its greater geometric complexity.

Overall, these results answer RQ3 by demonstrating that incorporating movable solar panels provides a substantial improvement in SRP model fidelity at an acceptable computational cost. Therefore, for spacecraft with large articulated structures, neglecting panel mobility can lead to significant orbit propagation errors, while the additional computational effort required to model these effects remains compatible with high-fidelity SRP simulations.

\section{Guidelines for Integrating High-Fidelity SRP Models into Orbit Propagation}

The experimental results indicate that the most appropriate SRP modeling strategy depends on the required balance between physical fidelity, propagation accuracy, and computational cost.

SPAD-based interpolation provides the most efficient solution for routine orbit propagation when the spacecraft attitude remains within the range represented by the precomputed dataset. Under these conditions, it reproduces the trajectory obtained from direct GPU ray tracing with only small deviations while reducing the computational cost by several orders of magnitude. Consequently, it is well suited for long-term propagation, mission analysis, and large-scale simulation campaigns. However, its accuracy depends on the coverage and resolution of the precomputed dataset, making direct SRP evaluation preferable when illumination conditions cannot be adequately represented through interpolation.

Direct GPU ray tracing should be employed when maximum SRP fidelity is required or when generating reference solutions and SPAD datasets. The proposed Vulkan implementation substantially reduces the computational cost of direct SRP evaluation compared with the previous OpenGL implementation, making high-fidelity ray tracing practical for more demanding propagation scenarios. Nevertheless, despite these improvements, direct evaluation remains considerably more expensive than interpolation-based approaches.

The results also show that the level of physical realism included in the spacecraft model should reflect the mission characteristics. For spacecraft with large articulated solar panels, neglecting panel motion can produce significant long-term propagation errors and should therefore be avoided. In contrast, for spacecraft with limited articulation or low SRP sensitivity, a simplified fixed-geometry model may provide an adequate compromise between accuracy and computational cost.

Overall, a practical workflow consists of generating high-fidelity reference data through direct GPU ray tracing, incorporating articulated structures whenever required, and using SPAD-based interpolation during routine orbit propagation whenever its validity conditions are satisfied.

These strategies are complementary rather than mutually exclusive, allowing the SRP modeling approach to be selected according to the accuracy, realism, and computational requirements of the mission.

\section{Conclusions and Future Work}

This work addressed three complementary objectives aimed at improving the practicality and physical realism of high-fidelity SRP simulations. First, interpolation-based SPAD datasets were evaluated as an efficient alternative to direct RayTrace GPU evaluation during orbit propagation. Second, a Vulkan-based implementation was developed to accelerate the computationally intensive ray-tracing process while preserving the numerical accuracy of the original model. Third, the simulation framework was extended to support movable solar panels, enabling a more realistic representation of spacecraft with articulated structures.

The evaluation of interpolation-based SRP propagation further showed that precomputed SPAD datasets can accurately reproduce the trajectories obtained with direct RayTrace GPU while reducing the computational cost of orbit propagation by several orders of magnitude. These results demonstrate that high-fidelity ray tracing can be effectively employed offline to generate SPAD datasets, enabling efficient long-term trajectory propagation with only a minor loss of accuracy.

The numerical validation confirmed that the Vulkan implementation preserves the fidelity of the original OpenGL-based RayTrace GPU model, with relative differences below $5\times10^{-4}$. From the computational perspective, performing both the ray tracing and the reduction on the GPU substantially reduces the execution time of direct SRP evaluations, achieving speed-ups of up to 9.4 for individual force computations and up to 15.2 for complete orbit propagation. These improvements demonstrate that direct high-fidelity SRP evaluation becomes practical even for complex spacecraft models, where the computational cost of the original implementation would otherwise be prohibitive.

The incorporation of movable solar panels increases the physical realism of the SRP model by accounting for the continuous variation of the illuminated spacecraft geometry. The results show that neglecting panel mobility can produce significant errors in both instantaneous SRP forces and long-term orbit propagation, particularly for spacecraft with large articulated solar arrays. Although modeling panel motion increases the ray-tracing workload, the resulting computational overhead remains moderate relative to the improvement in simulation fidelity.

Overall, the developments presented in this chapter significantly improve both the efficiency and realism of the proposed high-fidelity SRP simulation framework. The Vulkan implementation enables more efficient high-fidelity SRP simulations, while the movable-panel model extends the applicability of the framework to a wider range of spacecraft configurations without compromising the physical consistency of the SRP model.

Future work will focus on further optimizing the Vulkan implementation, particularly for lightweight spacecraft models where communication overhead remains significant. A key objective will be the integration of the Vulkan-based implementation with the movable solar panel model to combine both performance improvements and increased physical realism within a single simulation framework. Additional developments will include more advanced spacecraft articulation mechanisms and the incorporation of additional perturbation forces to further increase the fidelity of long-term orbit propagation.

Future investigations should also consider scenarios where SRP represents a more significant perturbation in the orbital dynamics, such as missions around small bodies or spacecraft with high area-to-mass ratios. In these cases, the differences between interpolation-based and direct ray-tracing approaches may become more pronounced, providing further insight into the trade-off between computational cost and SRP modeling fidelity.


\newpage
\begin{thebibliography}{9}


\bibitem{simulatorSRP}
L. Zardaín, A. Farrés, and A. Puig,
``Validation of the SRP tool: HiFi-SoRaP,''
in \textit{Proceedings of the 2023 AAS/AIAA Astrodynamics Specialist Conference},
Big Sky, MT, USA, Aug. 2023.

\bibitem{vulkan}
K. Muntal Freixas, 
\textit{Acceleration techniques for the computation of Solar Radiation Pressure}, Bachelor's Thesis, 2025. 

\bibitem{mobilePanels}
A. Mir Moncho,
\textit{Ray-Tracing Techniques for Solar Radiation Pressure Computation on Satellites with Mobile Pieces}, Bachelor's Thesis, 2025

\bibitem{rodriguez_rosseta}
T. Kato and J. C. van der Ha,
\textit{Precise modelling of solar and thermal accelerations on Rosetta},
Acta Astronautica,
vol. 72, pp. 165--177, 2012.

\bibitem{chang}
X. Chang, B. Männel, and H. Schuh,
\textit{An analysis of a priori and empirical solar radiation pressure models for GPS satellites},
Advances in Geosciences,
vol. 55,
pp. 33--45,
2021.

\bibitem{li}
F. Darugna, P. Steigenberger, O. Montenbruck, and S. Casotto,
\textit{Ray-tracing solar radiation pressure modeling for QZS-1},
Advances in Space Research,
vol.~62, no.~4, pp.~935--943, 2018.

\bibitem{zhang}
Z. Li, M. Ziebart, S. Bhattarai, D. Harrison, and S. Grey,
\textit{Fast solar radiation pressure modelling with ray tracing and multiple reflections},
\textit{Advances in Space Research},
vol.~61, no.~9, pp.~2352--2365, 2018.



\bibitem{Dormand1980AFO}
J.~R. Dormand and P.~J. Prince,
``A family of embedded Runge--Kutta formulae,''
\textit{Journal of Computational and Applied Mathematics},
vol.~6, pp.~19--26, 1980.

\bibitem{Hairer2006GNI}
E.~Hairer, C.~Lubich, and G.~Wanner,
\textit{Geometric Numerical Integration: Structure-Preserving Algorithms for Ordinary Differential Equations},
Springer, Berlin, Heidelberg, 2006.


\bibitem{Kenneally20}
P.~W. Kenneally and H.~Schaub,
``Fast spacecraft solar radiation pressure modeling by ray tracing on graphics processing unit,''
\textit{Advances in Space Research},
vol.~65, no.~8, pp.~1951--1964, 2020.

\bibitem{Vallado2001}
D.~A. Vallado,
\textit{Fundamentals of Astrodynamics and Applications},
Springer, Dordrecht, The Netherlands, 2001.


\bibitem{Montenbruck}
O.~Montenbruck and E.~Gill,
\textit{Satellite Orbits: Models, Methods, and Applications},
Springer, Berlin, Heidelberg, 2000.


\bibitem{Tichy2014FastFS}
J. Tichy, A. Brown, M. Demoret, B. Schilling, and D. Raleigh,
``Fast Finite Solar Radiation Pressure Model Integration Using OpenGL,''
in \textit{Proceedings of the 24th International Symposium on Space Flight Dynamics (ISSFD)},
Laurel, MD, USA, 2014.


\bibitem{HifiSorap}
L.~Zardaín, A.~Farrés, and A.~Puig,
``High-fidelity Modeling and Visualizing of Solar Radiation Pressure: A Framework for High-fidelity Analysis,''
in \textit{Proceedings of the 2020 AAS/AIAA Astrodynamics Specialist Conference},
Virtual Conference, Aug. 2020.


\bibitem{Kenneally16}
P.~W. Kenneally,
\textit{High Geometric Fidelity Solar Radiation Pressure Modeling via Graphics Processing Unit},
M.S. Thesis, University of Colorado Boulder, 2016.


\bibitem{kenneally2018brdf}
P.~W. Kenneally and H.~Schaub,
``Spacecraft Radiation Pressure Using Complex Bidirectional-Reflectance Distribution Functions on Graphics Processing Unit,''
in \textit{Proceedings of the 69th International Astronautical Congress},
Bremen, Germany, 2018.


\bibitem{VulkanProgrammingGuide}
G.~Sellers, J.~Kessenich, and D.~Shreiner,
\textit{Vulkan Programming Guide: The Official Guide to Learning Vulkan},
1st ed., Addison-Wesley Professional, 2017.

\end{thebibliography}
\end{document}